\documentclass[10pt]{article}
 \usepackage[T1]{fontenc}
 \usepackage{latexsym,amssymb,amsmath}
\usepackage{pstricks,lscape,rotating}
\usepackage{fancyhdr}

 \newcommand{\G}{\widetilde{G}}
 \newcommand{\Galg}{\mathbf{G}}
 
 \newcommand{\Balg}{\mathbf{B}}
 \newcommand{\redu}{\rho}
 
 \newcommand{\Ualg}{\mathbf{U}}
 \newcommand{\Talg}{\mathbf{H}}
 
 \newcommand{\Nalg}{\mathbf{N}}

 \newcommand{\RR}{\mathcal{R}}
 \newcommand{\uni}{\mathcal{U}}

 \newcommand{\Tr}{\operatorname{Tr}}
 \newcommand{\Gf}{{G}}
 \newcommand{\Su}{\operatorname{Sz}}

 \newcommand{\s}{\sigma}
 \newcommand{\Ind}{\operatorname{Ind}}
 \newcommand{\Res}{\operatorname{Res}}
 \newcommand{\Ker}{\operatorname{Ker}}
 
 \newcommand{\Cen}{\operatorname{C}}
 \newcommand{\Cl}{\operatorname{Cl}}
 \newcommand{\Sh}{\operatorname{Sh}}
 \newcommand{\St}{\operatorname{St}}
  
 \newcommand{\B}{\widetilde{B}}
 
 \newcommand{\Sz}{\widetilde{\Su}}
 
 \newcommand{\Uz}{\widetilde{U}_0}

 \newcommand{\semi}[1]{\rtimes\langle\,#1\,\rangle}
 \newcommand{\ps}[2]{\langle\,#1,#2\,\rangle}
 \newcommand{\cyc}[1]{\langle\,#1\,\rangle}

 \newcommand{\F}{\mathbb{F}}
 \newcommand{\K}{\operatorname{K}}
 \newcommand{\N}{\mathbb{N}}
 \newcommand{\Q}{\mathbb{Q}}
 \newcommand{\Z}{\mathbb{Z}}

 \newcommand{\Irr}{\operatorname{Irr}}

\title{\bf The Shintani descents of Suzuki Groups and consequences}

\author{Olivier \textsc{Brunat}\\
 brunat@igd.univ-lyon1.fr}

 \newtheorem{lemme}{Lemma}[section]
 
 \newtheorem{theoreme}{Theorem}[section]
 
 \newtheorem{proposition}{Proposition}[section]
 
 \newtheorem{remarque}{Remark}[section]

 \newtheorem{conjecture}{Conjecture}[section]
 \newenvironment{preuve}[1][]{\noindent {\bf{Proof }}{---\ }}{\hfill
 \nopagebreak $\Box$
 \newline
}

\newenvironment{changemargin}[2]{\begin{list}{}{%
\setlength{\topsep}{0pt}%
\setlength{\leftmargin}{0pt}%
\setlength{\rightmargin}{0pt}%
\setlength{\listparindent}{\parindent}%
\setlength{\itemindent}{\parindent}%
\setlength{\parsep}{0pt plus 1pt}%
\addtolength{\leftmargin}{#1}%
\addtolength{\rightmargin}{#2}%
}\item }{\end{list}}

\begin{document}
\maketitle

\begin{changemargin}{2cm}{2cm}
\begin{center}{\sc abstract}\end{center}
The main aim of this paper is to associate to every cuspidal unipotent character of the Suzuki group its root of unity and to give a possible definition of the Fourier matrix associated to the family of the cuspidal unipotent characters of this group. We compute to this end the Shintani descents of Suzuki groups and use results of Digne and Michel.  

\end{changemargin}

\section{Introduction}
Let $\Galg$ be a connected reductive group defined over the finite field with $q$ elements and let~$F:\Galg\rightarrow\Galg$ be a generalized Frobenius map. Let~$W$ be the Weyl group of~$\Galg$ which respect to an $F$-stable maximal torus of~$\Galg$ contained in an $F$-stable Borel subgroup of~$\Galg$. We denote by~$\Galg^F$ the finite group of fixed points under~$F$. For~$w\in W$ we have a corresponding generalized Deligne-Lusztig character~$R_w$ of~$\Galg^F$.
%(the definition of~$R_w$ is recall in~\S\ref{DL}).
We denote by~$\uni(\Galg^F)$ the set of unipotent characters of~$\Galg^F$, that is the irreducible constituent of the~$R_w$ for~$w\in W$. In~\cite{Lust} Lusztig attached to every~$\chi\in\uni(\Galg^F)$ a root of unity~$\omega_{\chi}$.
%(the definition of these roots is recall in \S\ref{root} of this paper).
In~\cite{Lusztig} he computes all such roots of unity corresponding to unipotent characters of
finite reductive groups, except for some pairs of complex conjugate
characters where a sign is missing. For example,
if~$\Galg^F=\Su(2^{2n+1})$ (where~$n$ is a non-negative integer)
is the Suzuki group with parameter~$2^{2n+1}$ then~$\Galg^F$ has
two cuspidal unipotent complex conjugate characters~$\mathcal{W}$
and~$\overline{\mathcal{W}}$ of degree~$2^n(2^{2n+1}-1)$;
%(these characters are precise in~??).
Lusztig's method only gives that~$\omega_{\mathcal{W}}=\frac{\sqrt{2}}{2}(-1\pm\sqrt{-1})$. %$$\{\omega_{\mathcal{W}},\omega_{\overline{\mathcal{W}}}\}=\{\frac{\sqrt{2}}{2}(-1-\sqrt{-1}),\frac{\sqrt{2}}{2}(-1+\sqrt{-1})\}.$$
On the other hand, using the {\it almost characters} of~$\Galg^F$, Lusztig has shown that the unipotent characters of~$\Galg^F$ can be distributed in families. %Let~$\mathcal{F}\subseteq\uni(\Galg^F)$ such a family.
In~\cite{Lust} Lusztig associated to most of the families a matrix,
the so-called Fourier matrix of the family. However the Suzuki and Ree groups do not have Fourier matrices in the Lusztig sense! On the other side,
Geck and Malle have axiomatized Fourier matrices and give
in~\cite{GM} candidates for these groups.

The aim of this paper is to compute the roots of unity attached to the cuspidal unipotent characters~$\mathcal{W}$ and~$\overline{\mathcal{W}}$ of the Suzuki group with parameter~$2^{2n+1}$ and to give a possible definition of the Fourier matrix associated to the family~$\{\mathcal{W},\overline{\mathcal{W}}\}$ with another approach as the one proposed by Geck and Malle. To this end, we compute {\it the Shintani descents} of the Suzuki group and use results of Digne-Michel~\cite{DM}.

Let~$n$ be a non-negative integer and let $\Galg$ be a simple
group of type $B_2$ defined over $\overline{\F}_2$.
Let~$F$ be the generalized Frobenius map such that~$\Galg^F$ is
the Suzuki group  with parameter~$2^{2n+1}$. The finite group~$\Galg^{F^2}=B_2(2^{2n+1})$ is
the finite "untwisted" group of type~$B_2$ with parameter~$2^{2n+1}$. 

This paper is organized as follows: in the
first section, we recall some definitions and generalities. In the
second section, we explicitly compute the character table of the
finite group~$B_2(2^{2n+1})\semi\s$, where $\s$ is the restriction
of $F$ to~$B_2(2^{2n+1})$. The main result of this part
is:\medskip

\begin{theoreme}\label{mainth}
Let~$n$ be a non-negative integer. We set~$q=2^{2n+1}$ and let~$\s$
the exceptional graph automorphism of~$B_2(q)$ such that its fixed
points subgroup is the Suzuki group with parameter~$q$. Then the
group~$B_2(q)\semi{\s}$ has~$(2q+6)$ irreducible extensions
of~$(q+3)$ irreducible $\s$-stable characters of~$B_2(q)$. The
values of these extensions are given in Table~\ref{CaractereB2}.
\end{theoreme}

In the last section, we compute the Shintani descents of the
Suzuki group with parameter~$2^{2n+1}$. We then obtain two consequences on the unipotent
characters of the Suzuki group: we first explicitly compute the root
of unity associated to their unipotent characters and secondly we
compute a Fourier matrix for the families of these groups.\\

I wish to express my hearty thanks to Meinolf Geck for leading me to this work and for valuable discussions.

\section{Generalities}
\subsection{Finite reductive groups}\label{DL}
%For generalities on reductive groups, we refer to~\cite{Carter2} and~\cite{Geck}.
Let $\Galg$ be a connected reductive group defined over the finite
field with $q=p^f$ elements. Let~$F$ be a generalized Frobenius
map over $\Galg$. We recall that the finite
subgroup~$\Galg^F=\{x\in\Galg\ |\ F(x)=x\}$ is a so-called {\it
finite reductive group}. Let~$\Talg$ be an~$F$-stable maximal torus
of~$\Galg$ contained in an $F$-stable Borel~$\Balg$ of~$\Galg$.
We set~$W=\operatorname{N}_{\Galg}(\Talg)/\Talg$ the Weyl group
of~$\Galg$. The map~$F$ induces an automorphism of~$W$ (also
denoted by~$F$ to simplify). We denote by~$\delta$ the order of this
automorphism.\\ We fix~$w\in W$ and we define the corresponding Deligne-Lusztig variety by:
$$X_w=\{x\Balg\ |\ x^{-1}F(x)\in\Balg w\Balg\}.$$ We recall that for every positive integer~$i$, we
can associate a $\overline{\Q}_{\ell}$-space to~$X_w$,
the $i$-th $\ell$-adic cohomology space with compact support
$H_c^i(X_w,\overline{\Q}_{\ell})$ over the algebraic
closure~$\overline{\Q}_{\ell}$ of the $\ell$-adic field (Here,~$\ell$ is a
prime not dividing~$q$). The
group~$\Galg^F$ acts on~$X_w$. This action induces a linear action
on $H_c^i(X_w,\overline{\Q}_{\ell})$. Thus these spaces
are~$\overline{\Q}_{\ell}\Galg^F$-modules. We define the generalized Deligne-Lusztig character by:
$$\forall g\in\Galg^F,\quad R_w(g)=\sum_{i\geq
0}(-1)^i\Tr(g,H_c^i(X_w,\overline{\Q}_{\ell})).$$ The set
of irreducible characters of~$\Galg^F$ is denoted
by~$\Irr(\Galg^F)$ and we denote by~$\cyc{\ ,\ }_{\Galg^F}$ the usual
scalar product on the space~$\Cen(\Galg^F)$ of the
$\overline{\Q}_{\ell}$-valued class functions
of~$\Galg^F$. We define the set~$\uni(\Galg^F)$ of unipotent
characters of~$\Galg^F$ by:
$$\uni(\Galg^F)=\{\chi\in\Irr(\Galg^F)\ |\ \exists w\in W,\ \cyc{R_w,\chi}_{\Galg^F}\neq 0\}.$$

\subsubsection{The root of a unipotent character}\label{root}
The group~$\cyc{F^{\delta}}$ acts on~$X_w$. This action induces a
linear endomorphism on~$H_c^i(X_w,\overline{\Q}_{\ell})$.
We also fix an eigenvalue~$\lambda$ of~$F^{\delta}$ on
$H_c^i(X_w,\overline{\Q}_{\ell})$ and we denote by~$F_{\lambda,i}$
its generalized eigenspace. The actions of~$\Galg^F$ and
of~$\cyc{F^{\delta}}$ on~$H_c^i(X_w,\overline{\Q}_{\ell})$
commute, thus~$F_{\lambda,i}$ is a $\overline{\Q}_{\ell}\Galg^F$-module. Moreover
the irreducible constituents which occur in the character associated
to this $\overline{\Q}_{\ell}\Galg^F$-module are unipotent characters of~$\Galg^F$.
Now let~$\chi\in\uni(\Galg^F)$. Then there exists $w\in W$,
$\lambda\in\overline{\Q}_{\ell}^{\times}$ and $i\in\N$ such that
$\chi$ occurs in the character associated to~$F_{\lambda,i}$.
Lusztig has shown that~$\lambda$, up to a power of~$q^{1/2}$, is a root of unity which depends only on~$\chi$ (denoted by~$\omega_{\chi}$). Thus there exists~$s\in\N$ such that~$\lambda=\omega_{\chi}q^{s/2}$~(see~\cite{DM}).

\subsubsection{Fourier matrices}\label{FourierMat}
%We not precisely define the Fourier matrices. We refer
We assume that~$\delta\neq 1$. We recall that $\cyc{F}$ acts on~$\Irr(W)$. More precisely if~$\rho\in\Irr(W)$, we define~$\rho^F$
by~$\rho^F(w)=\rho(F(w))$ for every~$w\in W$.
Let~$\rho\in\Irr(W)$ such that $\rho^F=\rho$, i.e., the inertial
group of~$\rho$ in~$W\semi{F}$ is $W\semi{F}$. It follows
that~$\rho$ has extensions to~$W\semi{F}$. Let~$\widetilde{\rho}$ be
such an extension; we define the {\it almost
character} associated to~$\widetilde{\rho}$ by:
$$\RR_{\widetilde{\rho}}=\frac{1}{|W|}\sum_{w\in
W}\widetilde{\rho}(w,F)R_w,$$ where the elements of~$W\semi F$ are denoted by~$(w,x)$ for every~$w\in W$ and~$x\in\cyc{F}$.  Let~$\chi,\,\chi'\in \uni(\Galg^F)$.
The characters~$\chi$ and~$\chi'$ are in the same family if and
only if there exists~$(\chi_i)_{i=1,\ldots, m}$,
where~$\chi_i\in\uni(\Galg^F)$ such that:
\begin{itemize}
\item[\textbullet]
We have $\chi_1=\chi$ and $\chi_m=\chi'$,
\item[\textbullet]
For every $2\leq i\leq m-1$, there exists an~$F$-stable character~$\rho_i\in\Irr(W)$ such that
$$\begin{array}{lll}\langle\,\chi_i,\RR_{\widetilde{\rho}_i}\,\rangle_{\Galg^F}
\neq 0&\textrm{and}&\langle\,\chi_{i+1},\RR_{\widetilde{\rho}_i}\,\rangle_{\Galg^F} \neq 0.
\end{array}$$
\end{itemize}
Let~$\mathcal{F}$ be a family of unipotent characters of~$\Galg^F$
obtained in this way. Except in the cases where~$\Galg^F$ is a Suzuki
group or a Ree group of type~$G_2$ or~$F_4$, Lusztig has shown that we
can associate a matrix~$M_{\mathcal{F}}$ to~$\mathcal{F}$. We
refer to~\cite{Lust} for details. In the case where~$\Galg^F$ is a
Suzuki group or a Ree group, Geck and Malle proposed in~\cite{GM}
candidates for Fourier matrices of these groups in agreement with
a general axiomatization of Fourier matrices that they developed.

\subsubsection{Shintani descents}\label{secshin}
We recall that the Lang map associated to a generalized Frobenius
map~$F$ is the map $L_F:\Galg\rightarrow\Galg, x\mapsto
x^{-1}F(x)$. Since~$\Galg^{F^{\delta}}$ is a finite $F$-stable
subgroup and~$F$ is an automorphism of the abstract
group~$\Galg$, it follows that~$F$ restricts to an automorphism
of~$\Galg^{F^{\delta}}$, also denoted by~$F$. The maps~$F$
and~$F^{\delta}$ have Lang's property, that is, their
associated Lang maps are surjective (see~\cite{Geck}
Th.~$4.1.12$). Using this fact, we can establish a
correspondence~$N_{F/F^2}$ between~$\Galg^F$
and~$\Galg^{F^{\delta}}\semi F$, the so-called {\it Shintani
correspondence}. More precisely, let~$g\in\Galg^F$; by the surjectivity of the Lang map~$L_{F^{\delta}}$, there exists~$x\in\Galg$ such
that~$g=L_{F^{\delta}}(x)$. Therefore we have~$(L_F(x^{-1}),F)\in\Galg^{F^{\delta}}\semi{F}$. This
correspondence induces a bijection between the conjugacy classes
of~$\Galg^F$ and the conjugacy classes
of~$\Galg^{F^{\delta}}\semi{F}$ which consist of elements of the
form~$(g,F)$, with~$g\in\Galg^{F^{\delta}}$. Moreover, we have:
\begin{equation}
\forall\,g\in
\Galg^{F},\quad|\Cen_{\Galg^{F^{\delta}}\semi{F}}\left(N_{F/F^{\delta}}(g)\right)|=\delta|\Cen_{\Galg^F}(g)|.
\label{eqcardinal}
\end{equation}
Using this correspondence, we can associate a class function of~$\Galg^F$ to every class function of~$\Galg^{F^{\delta}}\semi F$. Indeed let~$\psi\in\Cen(\Galg^{F^{\delta}}\semi F)$; we then define the {\it Shintani descent of~$\psi$} by~$\Sh_{F^{\delta}/F} \psi=\psi\circ N_{F/F^{\delta}}$.
We refer to~\cite{DM} for further details.

\subsubsection{The link between Shintani descents, Roots and Fourier matrices}
The set of the irreducible constituents of~$\Ind_{\Balg^{F^{\delta}}}^{\Galg^{F^{\delta}}}1_{ \Balg^{F^{\delta}}}$ is the so-called {\it principal series} of~$\Galg^{F^{\delta}}$. There exists a $1$-$1$ correspondence between the irreducible characters of~$W$ and the characters of the principal series of~$\Galg^{F^{\delta}}$ (see~\cite{Carter2}). Let~$\rho\in\Irr(W)$, then we denote by~$\chi_{\rho}$ its corresponding character. Similarly Malle has shown that there is a $1$-$1$ correspondence between~$\Irr(W\semi F)$ and the irreducible components of~$\Ind_{\Balg^{F^{\delta}}}^{\Galg^{F^{\delta}}\semi F}1_{ \Balg^{F^{\delta}}}$. %There exists a $1$-$1$ correspondence between the irreducible characters of~$W$ and the character of the principal serie of~$\Galg^{F^{\delta}}$.
We now assume that~$\rho=\rho^F$, therefore~$\rho$ has irreducible
extensions in~$W\semi{F}$. We fix such an
extension~$\widetilde{\rho}$ and we denote by~$\chi_{\widetilde{\rho}}$ the
corresponding character. Then~$\chi_{\widetilde{\rho}}$ is an
extension of~$\chi_{\rho}$ to~$\Galg^{F^{\delta}}\semi F$. We
have:
\begin{theoreme}(Digne-Michel \cite{DM})\label{thDM}
Let~$\rho\in\Irr(W)$ such that~$\rho^F=\rho$. Let~$\widetilde{\rho}\in\Irr(W\semi F)$ be an extension of~$\rho$. Then we have:
$$ \Sh_{F^{\delta}/F} \chi_{\widetilde{\rho}}=\sum_{V\in
\uni({\Galg^F})}\ps{\RR_{\widetilde{\rho}}}{V}_{\Galg^F}\omega_V V.$$
\end{theoreme}
\begin{remarque}
The Theorem is proved in~\cite{DM} in the case where~$F$ is a
Frobenius map. But the arguments are the same when~$F$ is a
generalized Frobenius map.
\end{remarque}

\noindent We now recall some conjectures of Digne and Michel (see~\cite{DM}):
\begin{conjecture}\label{conjDM}
Let~$\chi\in\uni(Galg^{F^{\delta}})$ such that~$\chi^F=\chi$. Let~$\widetilde{\chi}$ be an extension of~$\chi$ to~$\Galg^{F^{\delta}}\semi{F}$. Then:
\begin{enumerate}
\item The irreducible constituents of~$\Sh_{F^{\delta}/F}\widetilde{\chi}$ are unipotent characters of~$\Galg^F$ and lie in the same family~$\mathcal{F}$;
\item There exists a root of unity~$u$ such that
$$\pm u \Sh_{F^{\delta}/F}\widetilde{\chi}=\sum_{V\in\mathcal{F}} a_V\omega_V V.$$
In this case, the coefficients $a_V$ give (up to a sign) a row of the Fourier matrix associated to the family~$\mathcal{F}$.
\end{enumerate}
\end{conjecture}
\subsection{Suzuki groups}\label{B2}
Let~$\Galg$ be a simple group of type~$B_2$ defined
over~$\F_2$. The root system of~$\Galg$
is~$\Phi=\{-a,-b,-a-b,-2a-b,a,b,a+b,2a+b\},$ where~$\Pi=\{a,b\}$
is chosen as a fundamental  root system. We denote
by~$\Phi^+=\{a,b,a+b,2a+b\}$ the set of positive roots with respect to~$\Pi$. The Weyl group~$W$ of~$\Galg$ is the dihedral
group with~$8$ elements. We denote by~$x_r(t)$ ($r\in\Phi$,~$t\in
\overline{\F}_2$) the Chevalley generators. It is convenient to
identify~$\Galg$ with the symplectic group of dimension~$4$ over
the algebraic closure of~$\F_2$ defined by:
$$\operatorname{Sp}_4(\overline{\F}_2)=\{A\in\operatorname{M}_4(\overline{\F}_2)\
|\ ^t\!AJA=J\},\quad\textrm{where}\quad J=\begin{bmatrix}
  0 & 0 & 0 & 1 \\
  0& 0 & 1 & 0 \\
  0 & 1 & 0 & 0 \\
  1 & 0 & 0 & 0
\end{bmatrix}.$$
Representing matrices for the Chevalley generators are, for every~$t\in\overline{\F}_2$:
$$\begin{array}{cc} x_a(t)=
\begin{bmatrix}
  1 & t & 0 & 0 \\
  0& 1 & 0 & 0 \\
  0 & 0 & 1 & t \\
  0 & 0 & 0 & 1
\end{bmatrix}
,&x_b(t)=
\begin{bmatrix}
  1 & 0& 0 & 0 \\
  0& 1 & t & 0 \\
  0 & 0 & 1 & 0 \\
  0 & 0 & 0 & 1
\end{bmatrix},\\
\\
 x_{a+b}(t)=
\begin{bmatrix}
  1 & 0 & t & 0 \\
  0& 1 & 0 & t \\
  0 & 0 & 1 & 0 \\
  0 & 0 & 0 & 1
\end{bmatrix}
,&x_{2a+b}(t)=
\begin{bmatrix}
  1 & 0& 0 & t \\
  0& 1 & 0 & 0 \\
  0 & 0 & 1 & 0 \\
  0 & 0 & 0 & 1
\end{bmatrix}.
\end{array}
$$ Moreover, we have~$x_{-r}(t)=\ ^t\!x_r(t)$. We denote by~$\Talg=\{h(z_1,z_2)\ |\ z_1,z_2\in\overline{\F}_2\}$ the subgroup of diagonal matrices of~$\Galg$ and we set~$\Nalg=\Nalg_{\Galg}(\Talg)$. Then for the elements~$n_a$ and~$n_b$
of~$\Nalg$ we have the representing matrices
$$\begin{array}{cc} n_a=
\begin{bmatrix}
  0 & 1 & 0 & 0 \\
  1& 0 & 0 & 0 \\
  0 & 0 & 0 & 1 \\
  0 & 0 & 1 & 0
\end{bmatrix}
\quad\text{and}&n_b=
\begin{bmatrix}
  1 & 0& 0 & 0 \\
  0& 0 & 1 & 0 \\
  0 & 1& 0 & 0 \\
  0 & 0 & 0 & 1
\end{bmatrix}.
\end{array}$$
We recall the Chevalley relations of~$\Galg$; for every $u,\,
v,\, z_1,\, z_2 \in \overline{\F}_2,$ we have:
\begin{equation}\label{steinberg}
\begin{array}{lll}
x_a(u)x_b(v)&=&x_b(v)x_a(u)x_{a+b}(uv)x_{2a+b}(u^2v)\\
x_{a+b}(u)x_{2a+b}(v)&=&x_{2a+b}(v)x_{a+b}(u)\\
h(z_1,z_2)x_a(u)h(z_1,z_2)^{-1}&=&x_a(z_1uz_2^{-1})\\
h(z_1,z_2)x_b(u)h(z_1,z_2)^{-1}&=&x_b(z_2^2u)\\
h(z_1,z_2)x_{a+b}(u)h(z_1,z_2)^{-1}&=&x_{a+b}(z_1uz_2)\\
h(z_1,z_2)x_{2a+b}(u)h(z_1,z_2)^{-1}&=&x_{2a+b}(z_1^2u)\\
n_ah(z_1,z_2)n_a^{-1}&=&h(z_2,z_1)\\
n_bh(z_1,z_2)n_b^{-1}&=&h(z_1,z_2^{-1})\\
\end{array}
\end{equation}
Let~$n$ be a positive integer. We define~$F_{2^n}$ to be the Frobenius map with parameter~$2^n$ of~$\Galg$, hence it raises the coefficients of a matrix to their $2^n$-th powers. The
group~$\Galg$ has a graph endomorphism~$\alpha$ described in
Proposition~$12.3.3$ of~\cite{Carter}. It is given on generators
by:
$$
\begin{array}{lll}
 \alpha(x_a(t))&=&x_b(t^2)\\
 \alpha(x_b(t))&=&x_a(t)\\
 \alpha(x_{a+b}(t))&=&x_{2a+b}(t^2)\\
 \alpha(x_{2a+b}(t))&=&x_{a+b}(t)\\
 \alpha(h(z_1,z_2))&=&h(z_1\, z_2,z_1\, {z_2}^{-1})\\
 \alpha(n_a)&=&n_b\\
 \alpha(n_b)&=&n_a\\
\end{array}$$
We fix~$n$ a positive integer and we set~$\theta=2^n$ and~$q=2\theta^2$.
We define the map: $$F=F_{\theta}\circ \alpha.$$ Since~$F^2=F_{q}$, it follows that~$F$ is a generalized Frobenius map of~$\Galg$. The automorphism of~$W$ induced by~$F$ has order~$2$ (that is~$\delta=2$ with the preceding notations). The finite subgroup~$\Galg^F$ is the Suzuki group with parameter~$q$. Using~\cite{Ono}, this group is the same as the one studied in~\cite{Suzuki}. Moreover we have~$\Galg^{F^2}=\Galg^{F_q}=\operatorname{Sp}_{4}(\F_q)$. This is a finite "untwisted" group of type~$B_2$ with parameter~$q$. We denote by~$\Balg$ the Borel subgroup of the upper triangular matrices of~$\Galg$ and we set~$\Ualg$ to be the unipotent radical of~$\Balg$. The groups~$\Talg$,~$\Ualg$ and~$\Balg$ are~$F^2$-stable. We then define $H=\Talg^{F^2}$,~$U=\Ualg^{F^2}$ and~$B=\Balg^{F^2}$.
We denote by~$\s$ the restriction of~$F$ to~$\Galg^{F^2}$.
In the following we set~$\Galg^F=\Su(q)$ and $\Galg^{F^2}=\Gf$, and we remark that
$$\Gf^{\s}=\Su(q).$$
Now the aim is to obtain results about the unipotent characters of~$\Su(q)$ by using Shintani descents between~$\Gf\semi\s$ and~$\Su(q)$. Before we do this we must compute the irreducible characters of~$\Gf\semi\s$.

\section{The irreducible characters of~$B_2(q)\semi\s$.}
We use the notation of the preceding section. In this section we want
to compute the irreducible characters of the
extension~$\G=\Gf\semi\s$. This group is an extension by an automorphism
of~$\Gf$ of order~$2$. For generalities on character tables of
extensions by an automorphism of order~$2$ we refer
to~\cite{brunat1}\S 1. We recall some definitions and general
properties. The group~$\Gf$ is a normal subgroup of~$\G$. Thus a
conjugacy class of~$\G$ is either contained in~$\Gf$ or it has no
element in~$\Gf$. A class in the first case is called {\it an
inner class} and it is called {\it an outer class} in the second
case. A character~$\psi$ of~$\G$ is called {\it an outer
character} if there exists an outer element~$(g,\s)$ such
that~$\psi(g,\s)\neq 0$. We denote by~$\varepsilon$ the linear
character of~$\G$ with kernel~$\Gf$. Clifford theory shows that
the irreducible characters of~$\G$ can be parameterized by
the irreducible characters of~$\Gf$ as follows: let~$\chi\in\Irr(\Gf)$, then
either~$\chi^{\s}\neq\chi$ and~$\Ind_{\Gf}^{\G} \chi\in\Irr(\G)$,
or~$\chi^{\s}=\chi$ and~$\chi$ has two extensions in~$\G$ which
differ up to multiplication by~$\varepsilon$. Since the values of
extensions on~$(1,\s)$ are integers, in the case where this value is non-zero, we denote by~$\widetilde{\chi}$ the extension of~$\chi$ such
that~$\widetilde{\chi}(1,\s)>0$. 

\subsection{The outer classes of~$B_2(q)\semi\s$}

The Suzuki group~$\Su(q)$ has three maximal tori that are cyclic
groups:~$\cyc{\pi_0}$, $\cyc{\pi_1}$ and~$\cyc{\pi_2}$ of order~$(q-1)$,~$(q+2\theta+1)$ and~$(q-2\theta+1)$, respectively
(see~\cite{Suzuki}). We denote by~$E_0$ (resp.~$E_1$ and~$E_2$)
the set of non-zero classes modulo the equivalence relation~$\sim$
on~$\Z/(q-1)\Z$ (resp.~$\Z/(q+2\theta+1)\Z$
and~$\Z/(q-2\theta+1)\Z$) defined by $j\sim i\
\Longleftrightarrow\ j\equiv\pm i\mod(q-1)$ (resp.~$j\equiv \pm i,\pm qi\mod (q+2\theta+1)$
and~$j\equiv \pm i,\pm qi\mod(q-2\theta+1)$). We
put:~$$E=\{\pi_0^i,\,\pi_1^j,\,\pi_2^k\ |\ i\in E_0,\, j\in
E_1,\,k\in E_2\}.$$ The conjugacy classes of~$\Gf$ are recalled in
Table~\ref{ConjSp4} of the appendix. To simplify notation, we
denote~$x_r(1)$ by~$x_r$, where~$r\in\Phi$. We have:\medskip

\begin{theoreme}\label{tconj} Let~$n$ a non-negative integer. We
put~$\theta=2^n$ and~$q=2\theta^2$. Let~$\Gf=B_2(q)$ and~$\s$ the
exceptional automorphism of~$\Gf$ that defines~$\Su(q)$. Then the
group~$\G=B_2(q)\semi\s$ has~$(q+3)$ outer classes. The set
$$\{(1,\s),(x_a,\s),(x_{a+b},\s),(x_ax_{a+b},\s),(\pi,\s);\ \pi\in E\}$$
is a system of representatives of the outer classes of~$\G$.
Moreover, we have:
$$\begin{array}{lll}|\Cen_{\G}(1,\s)|&=&2q^2(q-1)(q^2+1),\\
|\Cen_{\G}(x_a,\s)|&=&4q,\\ |\Cen_{\G}(x_{a+b},\s)|&=&2q^2,\\
|\Cen_{\G}(x_ax_{a+b},\s)|&=&4q,\\
\Cen_{\G}(\pi_0^i,\s)&=&\cyc{\pi_0}\times\cyc{\s},\\
\Cen_{\G}(\pi_1^j,\s)&=&\cyc{\pi_1}\times\cyc{\s},\\
\Cen_{\G}(\pi_2^k,\s)&=&\cyc{\pi_2}\times\cyc{\s}.
\end{array}$$
\end{theoreme}\medskip

\begin{preuve}{} The Suzuki group with parameter~$q$ has~$(q+3)$
conjugacy classes (see~\cite{Suzuki}). Using the Shintani
correspondence, it follows that~$\G$ has~$(q+3)$ outer classes.
Using Table~\ref{ConjSp4}, we see that the elements of~$E$ are
not conjugate in~$\Gf$. Let~$x$ be in~$E$. Since~$\s(x)=x$, we
have:
$$\Cen_{\G}(x)=\Cen_{\Gf}(x)\semi{\s}.$$
Moreover,~$\Cen_{\Gf}(x)$ has odd order, hence~$|\Cen_{\G}(x)|_2=2$.
Thus~$\Cen_{\G}(x)$ has a unique class that consists of elements
of order~$2$. We choose~$(1,\s)$ as a representative of this class.
Using the $2$-Jordan decomposition of~$\G$ (see~\cite{brunat1}~Lemma~3.1),
it follows that the elemnts of~$\{(x,\s)\ |\ x\in E\}$ are not conjugate in~$\G$.
Moreover the group~$\Cen_{\Gf}(x)$ is abelian. Using
Lemma~3.2~in~\cite{brunat1} and the fact
that~$\Cen_{\Gf}(x)^{\s}\cap \cyc{(1,\s)}=\{1\}$, we deduce that:
$$\Cen_{\G}(x,\s)=\Cen_{\Gf}(x)^{\s}\times \cyc{(1,\s)}.$$
Furthermore~$\Cen_{\Gf}(x)^{\s}=\Cen_{\Su(q)}(x)$.
In~\cite{Suzuki} Prop.~16 it is proven that for every~$i\in
E_0$,~$j\in E_1$ and~$k\in E_2$, we
have~$\Cen_{\Su(q)}(\pi_0^i)=\cyc{\pi_0}$,~$\Cen_{\Su(q)}(\pi_1^j)=\cyc{\pi_1}$
and~$\Cen_{\Su(q)}(\pi_2^k)=\cyc{\pi_2}$. Thus we obtain~$(q-1)$
distinct outer classes of~$\G$.

Now we prove that~$(1,\s)$,~$(x_a,\s)$,~$(x_{a+b},\s)$
and~$(x_ax_{a+b},\s)$ are not conjugate in~$\G$. They are
of order~$2$,~$8$,~$4$ and~$8$ respectively. It then suffices to
prove that~$(x_a,\s)$ and~$(x_ax_{a+b},\s)$ are not conjugate
in~$\G$. Furthermore using the Bruhat decomposition of~$\Gf$, we
can show that two elements in~$(U,\s)$ are conjugate in~$\G$ if
and only if they are conjugate by an element of~$B$. Suppose there
exists~$b=uh\in B$ such that~$(b,1)(x_a,\s)=(x_ax_{a+b},\s)(b,1)$,
that is~$bx_a=x_ax_{a+b}\s(b).$ Let~$z_1,\,z_2\in\F_q$ such
that~$h=h(z_1,z_2)$ and we set~$z_0=z_1/z_2$. Then~$uhx_a=ux_a(z_0)h.$
Moreover:
$$\begin{array}{lll}
ux_a(z_0)&=&x_a(t_a)x_b(t_b)x_{a+b}(t_{a+b})x_{2a+b}(t_{2a+b})x_a(z_0)\\
&=&x_a(t_a+z_0)x_b(t_b)x_{a+b}(t_{a+b}+z_0
t_b)x_{2a+b}(t_{2a+b}+z_0^2 t_b),\\
x_ax_{a+b}\s(u)&=&x_a(1+t_b^{\theta})x_b(t_a^{2\theta})x_{a+b}(1+t_{2a+b}^{\theta}+t_a^{2\theta}t_b^{\theta})\\
&&x_{2a+b}(t_{a+b}^{2\theta}+t_a^{2\theta}t_b^{2\theta}).
\end{array}$$
By the uniqueness of the decomposition of the elements of~$B$, we deduce that $\s(h)=h$ and $$\left\{\begin{array}{lll} t_a+z_0&=&t_b^{\theta}+1\\
t_a^{2\theta}&=&t_b\\ z_0
t_b+t_{a+b}&=&t_a^{2\theta}t_b^{\theta}+t_{2a+b}^{\theta}+1\\
z_0^2 t_b+t_{2a+b}&=&t_a^{2\theta}t_b^{2\theta}+t_{a+b}^{2\theta}.
\end{array}\right.$$
Then we obtain~$z_0=1$ and~$z_1=z_2$. Furthermore since~$\s(h)=h$
we obtain~$z_1=z_2=1$. It follows that~$h=1$. Now we deduce from these
relations that~$t_b^{\theta}+t_b+1=0$ and
that~$t_b^{2\theta}+t_b+1=0$. Thus~$t_b^{\theta}=t_b^{2\theta}$
and~$t_b^{\theta}$ is a root of~$X^2+X$, that
is~$t_b^{\theta}\in\{0,1\}$. We then obtain a contradiction in the
relation~$t_b^{\theta}+t_b+1=0$. Therefore~$(x_a,\s)$
and~$(x_ax_{a+b},\s)$ are not conjugate in~$\G$.

We now compute the centralizer in~$\G$ of these elements. First we
remark that~$\Cen_{\G}(\s)=\Su(q)\times\cyc{\s}$.
Now let~$x\in\{x_a,x_{a+b},x_{2a+b}\}$. We therefore
have~$|\Cen_{\G}(x,\s)|=2|\{g\in \Gf\ |\ gx\s(g^{-1})=x\}|$. Using
the Bruhat decomposition and a similar calculation as above, we
prove that~$\{g\in \Gf\ |\ gx\s(g^{-1})=x\}\subseteq U$.
\begin{itemize}
\item The Chevalley relations give $$\left\{
\begin{array}{lll}
t_a&=&t_b^{\theta}\\t_{2a+b}^{\theta}+t_a^{2\theta}t_b^{\theta}&=&t_b+t_{a+b}\\
t_{a+b}^{2\theta}+t_a^{2\theta}t_b^{2\theta}&=&t_b+t_{2a+b}
\end{array}\right..$$
It follows that~$t_b\in\{0,1\}$. In both cases the number of
solutions is$$|\{(t_{2a+b},t_{a+b})\ |\
t_{2a+b}^{\theta}=t_{a+b}\}|=q.$$ Finally we have:
$|\Cen_{\G}(x_a,\s)|=|\Cen_{\G}(x_ax_{a+b},\s)|=4q.$
\item We have $$\{g\in \Gf\ |\ gx_{a+b}\s(g^{-1})=x_{a+b}\}=\{u\in U\ |\
\s(u)=u\}=U\cap \Su(q).$$ Furthermore~$|U\cap \Su(q)|=q^2$. We thus
deduce that~$|\Cen_{\G}(x_{a+b},\s)|=2q^2$.
\end{itemize}

\end{preuve}

\subsection{The~$\s$-stable characters of~$B_2(q)$}

In Table~\ref{CarSp4} of the appendix we recall the values of
the irreducible characters of~$B_2(q)$ that we need in this work.
We have:
\begin{proposition}\label{carstableB2}
The group~$B_2(q)$ has~$(q+3)$ $\s$-stable irreducible characters:
\begin{itemize}
\item[\textbullet] The $\s$-stable unipotent characters~$1_{\Gf}$,~$\theta_1$,~$\theta_4$
and~$\theta_5$ of degree
respectively~$1$,~$\frac{1}{2}q(q+1)^2$,~$q^4$
and~$\frac{1}{2}q(q-1)^2$.
\item[\textbullet] The $\frac{1}{2}(q-2)$ characters~$\chi_1(i, (2\theta-1)i)$,~$i\in E_0$ of degree~$(q+1)^2\,
 (q^2+1)$. We denote these characters by~$\chi_{\pi_0}(i)$.
\item[\textbullet] The $\frac{1}{4}(q+2\theta)$ characters~$\chi_5((q-2\theta+1)j)$,~$j\in E_1$
of degree~$(q^2-1)^2$, denoted by~$\chi_{{\pi}_1}(j)$.
\item[\textbullet] The $\frac{1}{4}(q-2\theta)$ characters~$\chi_5((q+2\theta+1)k)$,~$k\in E_2$ of degree~$(q^2-1)^2$, denoted
by~$\chi_{\pi_2}(k)$.

\end{itemize}
\end{proposition}

\subsection{Irreducible characters obtained by induction from $B\semi{\s}$} The group~$B$ is~$\s$-stable,
thus~$\B=B\semi\s\subseteq \G$. We now induce some characters
of~$\B$ to~$\G$ which permit to obtain the outer values
of~$\widetilde{\chi}_{\pi_0}(i)$ ($i\in E_0$) and
of~$\widetilde{\theta}_4$. Let~$\gamma_0$ the primitive $(q-1)$-th
root of unity given in Table~\ref{CarSp4}.
We define the primitive
$(q-1)$-th root of unity~$\varepsilon_0=\gamma_0^{(4-4\theta)}$ and
we set~$\varepsilon_0^i(\pi_0^l)=\varepsilon_0^{il}+\varepsilon_0^{-il}$.

\begin{proposition}\label{chii}
We have:
\label{valeurchi}$$\renewcommand{\arraystretch}{1.4}
\begin{array}{c|c|c|c|c|c|c|c}
 &(1,\s)&(x_a,\s)&(x_{a+b},\s)&(x_ax_{a+b},\s)&(\pi_0^l,\s)&(\pi_1^m,\s)&(\pi_2^r,\s)\\ \hline
\widetilde{\chi}_{\pi_0}(i)&q^2+1&1&1&1&\varepsilon_0^{i}(\pi_0^l)&0&0\\
\widetilde{\theta}_4&q^2&0&0&0&1&-1&-1
\end{array}$$
\end{proposition}

\begin{preuve}{}
We have~$U\lhd\B$ and we denote by~$\pi_U:\B\rightarrow\B/U$ the
canonical map. Let~$\phi\in\Irr(\B/U)$; then~$\phi\circ\pi_U$ is
an irreducible character of~$\B$. We have~$\B/U\simeq
H\rtimes\langle\s\rangle$. We now construct outer characters
of~$H\rtimes\langle\s\rangle$. Since~$H \simeq
\F_q^{\times}\times\F_q^{\times}$, it follows that the irreducible
characters of~$H$ are:
 $$\forall\ 1\leq k,\,l\leq q-1 \hspace{10pt}
\phi_{k,l}(\gamma^i,\gamma^j)=\gamma_0^{ik+jl}.$$ The $\s$-stable
irreducible characters of~$H$ are~$\phi_{i,(2\theta-1)i}$ ($
i\in\{1,\ldots,q-1\}$). Using~\cite{brunat1} Lemma~3.5, we
construct~$(q-2)$ linear characters of~$H\rtimes\langle\s\rangle$
defined by~$\varphi_i(h,x)=\phi_{i,(2\theta-1)i}(h)$,
where~$x\in\{1,\s\}$. We write~$\phi_{i,\B}=\varphi_i\circ\pi_U$; we
thus obtain~$(q-2)$ linear characters of~$\B$. Moreover~$\B$
has~$(q+2)$ outer classes, which are~$(\pi_0^l,\s)$
($l\in\{1,(q-2)\}$) and~$(1,\s)$,~$(x_a,\s)$,~$(x_{a+b},\s)$
and~$(x_ax_{a+b},\s)$. It is then easy to compute the values
of~$\phi_{i,\B}$. Moreover using the Mackey formula, it follows
that~$\Ind_{\B}^{\G}\phi_{i,\B}$ is an irreducible extension
of~${\chi}_{\pi_0}(i)$. To obtain the outer values of these
characters, we will induce~$\phi_{i,\B}$ from~$\B$ to~$\G$. To
this end, we give the corresponding induction formula. Except
for~$(1,\s)$, the centralizers of the outer elements of~$\B$ are
the same as their centralizers in~$\G$. We
have~$\Cen_{\B}(1,\s)=(\Su(q)\cap B,\s)=B^{\s}\times\cyc{\s}$ of order~$2q^2(q-1)$ and
$$\Cl_{\G}(\pi_0,\s)\cap\B=\Cl_{\B}(\pi_0,\s)\cup\Cl_{\B}(\pi_0^{-1},\s).$$
This permits to obtain the induction formula from~$\B$ to~$\G$
given in the following table:
\label{valeurinduitB}\small$$\renewcommand{\arraystretch}{1.4}
\begin{array}{c|c|c|c|c|c}
 &(1,\s)&(x_a,\s)&(x_{a+b},\s)&(x_ax_{a+b},\s)&(\pi_0,\s)\\ \hline
\Ind_{\B}^{\G}\phi&(q^2+1)\phi(1,\s)&\phi(x_a,\s)&\phi(x_{a+b},\s)&\phi(x_ax_{a+b},\s)
&\phi(\pi_0,\s)+\phi(\pi_0^{-1},\s)
\end{array}$$\normalsize
Using this formula, we compute the values
of~$\widetilde{\chi}_{\pi_0}(i)$.

\noindent We have:
$$\begin{array}{lllllll} \langle \Ind_B^{\Gf} 1_B
,1_G\rangle_{\Gf}&=&1&&\langle \Ind_B^{\Gf} 1_B, \theta_1\rangle_{\Gf}&=&2\\
\langle \Ind_B^{\Gf} 1_B, \theta_4\rangle_{\Gf}&=&1&&\langle
\Ind_B^{\Gf} 1_B, \theta_5\rangle_{\Gf}&=&0\\ \langle \Ind_B^{\Gf}
1_B, \chi\rangle_{\Gf}&=&0 &&\forall\, \chi=\chi^\s
\hspace{10pt}\text{such} &\text{that }& \chi \neq 1,\,\theta_4,\,
\theta_1,\, \theta_5.
\end{array}$$
We have~$\Res_{\Gf}^{\G}(\Ind_{\B}^{\G}1_{\B})=\Ind_B^G 1_B$ and write~$\xi$ for  the constituent
of~$\Ind_{\B}^{\G}1_{\B}$ which is an extension of~$\theta_1$.
Since~$\cyc{\Ind_B^{\Gf} 1_B, \theta_1}_{\Gf}=2$, it follows
that~$\xi$ occurs in~$\Ind_{\B}^{\G}1_{\B}$ with
multiplicity~$1$ or~$2$. If this multiplicity is~$2$, then
$$\langle \Ind_{\B}^{\G}1_{\B}, \Ind_{\B}^{\G}1_{\B}\rangle_{\G}
\geq 6.$$ But a direct calculation gives~$\langle \Ind_{\B}^{\G}1_{\B},
\Ind_{\B}^{\G}1_{\B}\rangle_{\G}=5.$ Thus it follows
that~$\xi+\varepsilon\xi$ is a
constituent of~$\Ind_{\B}^{\G}1_{\B}$. Decomposing the
character~$\Ind_{\B}^{\G}(1_{\B})-1_{\G} -
\xi-\varepsilon\xi$, we find an
irreducible extension of~$\theta_4$. We use the preceding
induction formula to compute its values.

\end{preuve}

\subsection{Irreducible characters obtained by induction from\\ $\Su(q)\times\cyc{\s}$}
Let~$\rho_0=x_ax_{a+b}\in\Su(q)$ and~$\sigma_0=x_{a+b}x_{2a+b}\in\Su(q)$. Let~$\tau_0$ be the complex primitive root of order~$q^2+1$ which appears in the character table of~$\Gf$. We recall that~$\{1,\sigma_0,\rho_0,\rho_0^{-1},\pi_0^i;i\in E_0,\pi_1^j;j\in E_1,\pi_2^k;k\in E_2\}$ is a system of representatives of the classes of~$\Su(q)$ (see~\cite{Suzuki}). We put~$\varepsilon_1=\tau_0^{(q-2\theta+1)^2}$ and~$\quad\varepsilon_2=\tau_0^{(q+2\theta+1)^2}$.
The character Table of~$\Su(q)$ is computed in~\cite{Suzuki} and reprinted in the appendix for the convenience of the reader. 
Since~$\Su(q)$ is the subgroup of fixed points under~$\s$, it follows that~$\Sz(q)=\Su(q)\times\cyc{\s}\subseteq\G$. 
Thus the classes and the character table of~$\Sz(q)$ are directly obtained using the classes and the character table of~$\Su(q)$. 
We give the induction formula from~$\Sz(q)$ to~$\G$ in Table~\ref{tabindSz}.
\begin{table}
$$\begin{array}{cc}
\renewcommand{\arraystretch}{1.5}\begin{array}{c|l}
 &\Ind_{\Su(q)}^{G}\phi\\\hline
A_1&q^2(q+1)(q^2-1)\phi(1)\\ A_{32}&q^2\phi(\s_0) \\
A_{42}&q(\phi(\rho_0)+\phi(\rho_0^{-1}))\\ \pi_0&(q-1)\phi(\pi_0)\\
\pi_1&(q-r+1)\phi(\pi_1) \\ \pi_2&(q+r+1)\phi(\pi_2)
\end{array}$$
&
$$\renewcommand{\arraystretch}{1.5}\begin{array}{c|l}
 &\Ind_{\Sz(q)}^{\G}\phi\\ \hline
(\pi_0,\s)&\phi(\pi_0,\s)\\
(1,\s)&\phi(1,\s)+(q-1)(q^2+1)\phi(\s_0,\s)\\ (x_a,\s)&0\\
(x_{a+b},\s)&\frac{q}{2}(\phi(\rho_0,\s)+\phi(\rho_0^{-1},\s))\\
(x_ax_{a+b},\s)&0\\ (\pi_1,\s)&\phi(\pi_1,\s)\\
(\pi_2,\s)&\phi(\pi_2,\s)
\end{array}
\end{array}$$
\caption{Induction formula from~$\Sz(q)$ to~$\G$.\label{tabindSz}}
\end{table}
Let~$\widetilde{\phi}\in\Irr(\Sz(q))$; to simplify we denote by the same symbol its induced character of~$\G$.
\begin{proposition} We set~$\varepsilon_0^l(\pi_0^i)=\varepsilon_0^{il}+\varepsilon_0^{-il}$,~$ \varepsilon_1^l(\pi_1^j)=\varepsilon_1^{jl}+\varepsilon_1^{-jl}+
\varepsilon_1^{jlq}+\varepsilon_1^{-jlq}$ and~$\varepsilon_2^l(\pi_2^k)=\varepsilon_2^{kl}+\varepsilon_2^{-kl}+
\varepsilon_2^{klq}+\varepsilon_2^{-klq}$.
In Tables~\ref{TabIndgf} and~\ref{TabIndgftild}, we give the values of the induced characters from~$\Sz(q)$ to~$\G$.
\end{proposition}
\begin{table}
\small$$\renewcommand{\arraystretch}{1.5}\begin{array}{c|l|l|l}
 &A_1&A_{32}&A_{42}\\ \hline
\Ind_{\Su(q)}^{\Gf}1_{\Su(q)}&q^2(q+1)(q^2-1)&q^2&2q\\
\Ind_{\Su(q)}^{\Gf}\St&q^4(q+1)(q^2-1)&0&0\\
\Ind_{\Su(q)}^{\Gf}X_i&q^2(q^2+1)(q+1)(q^2-1)&q^2&2q\\
\Ind_{\Su(q)}^{\Gf}Y_j&q^2(q-\theta+1)(q^2-1)^2&q^2(\theta-1)&-2q\\
\Ind_{\Su(q)}^{\Gf}Z_k&q^2(q+\theta+1)(q^2-1)^2&-q^2(\theta+1)&-2q\\
\Ind_{\Su(q)}^{\Gf}\mathcal{W}&\frac{1}{2}q^2\theta(q^2-1)^2&-q^2\theta&0\\ \multicolumn{4}{c}{}\\
&\pi_0&\pi_1&\pi_2\\ \hline \Ind_{\Su(q)}^{\Gf}1&q-1&q-\theta+1&q+\theta+1\\
\Ind_{\Su(q)}^{\Gf}\St&q-1&-(q-\theta+1)&-(q+\theta+1)\\
\Ind_{\Su(q)}^{\Gf}X_i&(q-1)(\varepsilon_0^i\pi_0)&0&0\\
\Ind_{\Su(q)}^{\Gf}Y_j&0&-(q-\theta+1)\varepsilon_1^j(\pi_1)&0\\
\Ind_{\Su(q)}^{\Gf}Z_k&0&0&-(q+\theta+1)\varepsilon_2^k(\pi_2)\\
\Ind_{\Su(q)}^{\Gf}\mathcal{W}&0&q-\theta+1&-(q+\theta+1)\\
\end{array}$$\\
\caption{Values of the induced character of~$\Su(q)$.\label{TabIndgf}}
\end{table}
\begin{table}
$$\renewcommand{\arraystretch}{1.5}\begin{array}{c|l|l|l|l|l}
 &(1,\s)&(x_{a+b},\s)&(\pi_0,\s)&(\pi_1,\s)&(\pi_2,\s)\\ \hline
\widetilde{1}_{\Su(q)}&q(q^2-q+1)&q&1&1&1\\ \widetilde{\St}&q^2&0&1&-1&-1\\
\widetilde{X_i}&q(q^2+1)&q&\varepsilon_0^i(\pi_0)&0&0\\
\widetilde{Y_j}&q-1)(q+q^2(\theta-1))&-q&0&-\varepsilon_1^j(\pi_1)&0\\
\widetilde{Z_k}&(q-1)(q-q^2(\theta+1))&-q&0&0&-\varepsilon_2^k(\pi_2)\\
\widetilde{\mathcal{W}}&-\frac{1}{2}q^2\theta(q-1)&0&0&1&-1\\
\end{array}$$
\caption{Outer values of the induced characters of~$\Sz(q)$\label{TabIndgftild}}
\end{table}
Let~$\psi$ be a generalized character of~$\G$. We recall that we can associate to~$\psi$ its~$\s$-reduction~$\rho(\psi)$, which is a character of~$\G$ such that~$\rho(\psi)(g,\s)=\psi(g,\s)$ (for every~$g\in\Gf$) and the irreducible constituents of~$\rho(\psi)$ are outer characters of~$\G$. We refer to~\cite{brunat1}~\S3.2.1 for details.
We now define:
$$\begin{array}{lcl}
X_0&=&\displaystyle{\redu\left(\widetilde{1}_{\Su(q)}-1_{\G}-\widetilde{\theta_4}-\sum_{i\in E_0}\widetilde{\chi}_{\pi_0}(i)\right),}\\ 
W_0&=&\displaystyle{\redu\left(\varepsilon\widetilde{W}-\theta(\widetilde{1}_{\Su(q)}-1_{\G})\right)}.
\end{array}$$

\begin{proposition}\label{propX1} For every~$k\in E_1$ (resp.~$k\in E_2$), there exists an extension~$\psi_k$ (resp.~$\psi'_k$)
of~$\chi_{\pi_1}(k)$ (resp.~$\chi_{\pi_2}(k)$) such that:
%On a\,\footnote{Les caractères~$\widetilde{\chi}_{A_1}(k)$
%et~$\widetilde{\chi}_{A_2}(k)$ désignent {\it a priori} une
%extension de~$\chi_{A_1}(k)$ et de~$\chi_{A_2}(k)$ respectivement.
%On va voir par la suite (cf.~la preuve de la
%proposition~\ref{theta1}) que les constituants qui apparaissent
%sont exactement les extensions données par la convention que l'on
%s'est fixé dans~\S\ref{genetable}. On adopte volontairement cet
%abus de notation, dans un souci de clareté.}:
$$\begin{array}{lcl}X_0&=&\displaystyle{\sum_{k\in E_1}\psi_k+\sum_{k\in E_2}\psi'_k,}\\&\\
W_0&=&\displaystyle{\sum_{k\in E_1}\psi_k-\sum_{k\in
E_2}\psi'_k}.\end{array}$$ We give the values of~$X_0$
and~$W_0$ in the following table: \small
$$\begin{array}{c|c|c|c|c|c|c|c|c}
 &(1,\s)&(x_a,\s)&(x_{a+b},\s)&(x_ax_{a+b},\s)&(\pi_0,\s)&(\pi_1,\s)&(\pi_2,\s)\\ \hline
 &&&&&&&\\
X_0&q(q-1)^2/2&-q/2&q/2&-q/2&0&1&1\\&&&&&&&\\
W_0&-\theta(q-1)&\theta&-\theta(q-1)&\theta&0&-1&1\\
\end{array}$$
\normalsize
\end{proposition}

\begin{preuve}{}
By computing the scalar products of~$\widetilde{1}_{\Su(q)}$ with the irreducible characters obtained in Prop.~\ref{chii}, we deduce that~$\widetilde{1}_{\Su(q)}-1_{\G}-\widetilde{\theta}_4-\sum_{i\in E_0}\widetilde{\chi}_{\pi_0}(i)$ is a character. Furthermore, we have~$$\begin{array}{l}\langle\Ind_{\Su(q)}^{\Gf}(1_{\Su(q)}) , \theta_1\rangle_{\Gf}=\langle\Ind_{\Su(q)}^{\Gf} (1_{\Su(q)}),\theta_5 \rangle_{\Gf}=0\\
\langle\Ind_{\Su(q)}^{\Gf}(1_{\Su(q)}) , \chi_{\pi_1}(k)\rangle_{\Gf}=\langle\Ind_{\Su(q)}^{\Gf} (1_{\Su(q)}) ,
\chi_{\pi_2}(k)\rangle_{\Gf}=1.\end{array}$$ This proves that~$X_0$ is exactly the summand of one extension of~$\chi_{\pi_1}(k)$ (denoted by~$\psi_k$) and one extension of~$\chi_{\pi_2}(k)$ (denoted by~$\psi'_k$). Now we compute the scalar products of~$\widetilde{\mathcal{W}}$ with the known irreducible characters of~$\G$. Thus~$\widetilde{\mathcal{W}}-\theta(\widetilde{\theta}_4\varepsilon+\sum\widetilde{\chi}_{\pi_0}(i)\varepsilon)$ is a character of~$\G$ and we denote by~$\Theta$ its~$\s$-reduction. We compute that~$\langle\Res_{\Gf}^{\G} \widetilde{\mathcal{W}},\theta_1
\rangle_{\Gf}=\langle\Res_{\Gf}^{\G}
\widetilde{\mathcal{W}},\theta_5 \rangle_{\Gf}=0$, $
\langle\Res_{\Gf}^{\G}\widetilde{\mathcal{W}} ,
\chi_{\pi_1}(k)\rangle_{\Gf}=\theta-1$ and~$
\langle\Res_{\Gf}^{\G}\widetilde{\mathcal{W}},
\chi_{\pi_2}(k)\rangle_{\Gf}=\theta+1$.
Since~$X_0$ and~$\widetilde{\mathcal{W}}$ have no common constituents (because~$\cyc{X_0,\widetilde{\mathcal{W}}}_{\G}=0$) we deduce that~$\Theta=(\theta-1)\sum_{k\in E_1}\psi_k\varepsilon+
(\theta+1)\sum_{k\in E_2}\psi'_k\varepsilon$. We remark that~$W_0=\Theta\varepsilon-\theta X_0$ and deduce that $$W_0=\sum\limits_{k\in E_2}\psi'_k-\sum\limits_{k\in E_1}\psi_k.$$

\end{preuve}
For every~$j\in E_1$ and~$k\in E_2$, we put~$
\varphi_j=\redu(\widetilde{Y}_j-\widetilde{Y}_1)$ and~$\vartheta_k=\redu(\widetilde{Z}_k-\widetilde{Z}_1)$.

\begin{proposition}\label{theta1}
Using the preceding notation we have:$$\varepsilon\widetilde{\chi}_{\pi_1}(1)=\frac{4}{q+2\theta}\left(\sum_{j\in
E_1}\varphi_j +\frac{1}{2}(X_0-W_0)\right),$$
$$\varepsilon\widetilde{\chi}_{\pi_2}(1)=\frac{4}{q-2\theta}\left(\sum_{k\in E_2}\vartheta_k
+\frac{1}{2}(X_0+W_0)\right).$$
Consequently we deduce that:
$$\begin{array}{c|c|c|c|c|c|c}
 &(1,\s)&(x_a,\s)&(x_{a+b},\s)&(x_ax_{a+b},\s)&(\pi_1,\s)&(\pi_2,\s)\\ \hline
 &&&&&&\\
\widetilde{\chi}_{\pi_1}(j)&(q-1)(q-2\theta+1)&-1&2\theta-1&-1&-\varepsilon_1^j(\pi_1)&0\\&&&&&&\\
\widetilde{\chi}_{\pi_2}(k)&(q-1)(q+2\theta+1)&-1&-2\theta-1&-1&0&-\varepsilon_2^k(\pi_2)\\
\end{array}$$
\end{proposition}

\begin{preuve}{}
We have~$ \langle\, \Res_{\Gf}^{\G}(\varphi_j),
\chi_{\pi_1}(j)\,\rangle_{\Gf}=-1$ and~$\langle
\,\Res_{\Gf}^{\G}(\varphi_j), \chi_{\pi_1}(1)\,\rangle_{\Gf}=1$. Moreover for every other $\s$-stable character~$\chi$ of~$\Gf$ we have~$\langle\, \Res_{\Gf}^{\G}(\widetilde{\varphi_j}),
\chi\,\rangle_{\Gf}=0$. This yields~$\Res_{\Gf}^{\G}(\varphi_j)=\chi_{\pi_1}(1)-\chi_{\pi_1}(j)$. We similarly prove that~$\Res_{\Gf}^{\G}(\vartheta_k)=\chi_{\pi_2}(1)-\chi_{\pi_2}(k)$.
We denote by~$\theta_j$ (resp.~$\theta_{1,j}$) the extension of~$\chi_{\pi_1}(j)$ (resp.~$\chi_{\pi_1}(1)$), which is an irreducible component of~$\varphi_j$. In particular we have~$\varphi_j=\theta_{1,j} -\theta_j$. First we prove that~$\theta_{1,j}$ is independent of~$j$. Indeed, we immediately compute that~$\langle\, \varphi_j - \varphi_2,\,\varphi_j -
\varphi_2\,\rangle_{\G}=2$ (where~$j\geq 3$). Furthermore we have~$
\varphi_j-\varphi_2=\theta_2-\theta_j+\theta_{1,j}-\theta_{1,2}$. Moreover since~$\Res_{\Gf}^{\G}(\varphi_j-\varphi_2)=\chi_{\pi_1}(2)-\chi_{\pi_1}(j),$ the characters~$\theta_2$ and~$\theta_j$ are constituents of~$\varphi_j-\varphi_2$. Since~~$\varphi_j-\varphi_2$ has two constituents, it follows that~~$\theta_{1,j}-\theta_{1,2}=0$, i.e., $$\forall\,j\geq 2\quad
\theta_{1,j}=\theta_{1,2}.$$ We denote by~$\theta_1$ this common constituent. We compute that~$\langle \varphi_i,X_0\rangle_{\G}=0$.
Also, if~$\theta_1=\psi_1$ then for every~$j\geq 2$, we have~$\theta_j=\psi_j$.
Indeed, if this is not the case we have~$\langle \varphi_i,X_0\rangle_{\G}=1
\neq 0$. With a similar argument we prove that if~$\theta_1=\varepsilon\psi_1$ then for every~$j\geq 2$, we obtain~$\theta_j=\varepsilon\psi_j$. Insummary we have:
$$\left\{
\begin{array}{lllll}
\theta_j&=&\psi_j&&\forall\ j\in E_1\\  &\text{or}&&&\\
\theta_j&=&\psi_j\varepsilon&&\forall\ j\in E_1
\end{array}\right..$$
Moreover, $$\frac{1}{4}(q+2\theta)\theta_1=\frac{1}{4}(q+2\theta)\theta_1
-\sum_{j\in E_1}\theta_j+\sum_{i\in E_1}\theta_j=\sum_{j\neq
1}\varphi_j +\sum_{j\in E_1}\theta_j.$$ Thus
$$\theta_1=\frac{4}{q+2\theta}\left(\sum_{j\neq 1}\varphi_j+\sum_{j\in E_1}\theta_j\right).$$
We immediately deduce from Proposition~\ref{propX1} that, $$\sum_{j\in
E_1}\psi_j=\frac{1}{2}(X_0-W_0).$$ Hence we have,
$$\left\{\begin{array}{lll}
\displaystyle{\sum_{j\in E_1}\theta_j}&=&\displaystyle{\frac{1}{2}(X_0-W_0)}\\
&\text{or}&\\
\displaystyle{\sum_{j\in
E_1}\theta_j}&=&\displaystyle{\frac{1}{2}(X_0-W_0)\varepsilon}
\end{array}\right..$$
We set $$\left\{\begin{array}{lll}
f_1&=&\displaystyle{\frac{4}{q+2\theta}\left(\sum_{j\neq
1}\varphi_j+\frac{1}{2}(X_0-W_0)\right),}\\&&\\
f_2&=&\displaystyle{\frac{4}{q+2\theta}\left(\sum_{j\neq
1}\varphi_j+\frac{1}{2}(X_0-W_0)\varepsilon\right).}\end{array}\right.$$
Either~$f_1$ or~$f_2$ is an irreducible character of~$\G$. But~$f_1(\pi_1,\s)$ is not an algebraic integer. Thus~$f_1$ is not a character. It follows that~$f_2=\varepsilon\psi_1$.
We similarly prove that
$$\varepsilon\psi'_1=\frac{4}{q-2\theta}\left(\sum_{k\neq 1}\vartheta_k
+\frac{1}{2}(X_0+W_0)\varepsilon\right).$$
Since we have~$\psi_j= \varepsilon\phi_j$, we immediately deduce the values of~$\psi_j$
using the relation~$\psi_j=\psi_1-\varepsilon\varphi_j$. Similarly we compute the values of~$\psi'_k$ using the relation~$\psi'_k=\psi'_1-\varepsilon\vartheta_k$.

\end{preuve}

The induced characters of~$\Sz(q)$ are not sufficient to obtain the outer values of the extensions of~$\widetilde{\theta}_1$ and~$\widetilde{\theta}_5$. However decomposing~$\widetilde{Y_i}$ we obtain:
 
\begin{lemme}\label{sommes}
The outer values of~$\widetilde{\theta}_1+\widetilde{\theta}_5$ are:\small
$$\begin{array}{c|c|c|c|c|c|c|c}
 &(1,\s)&(x_a,\s)&(x_{a+b},\s)&(x_ax_{a+b},\s)&(\pi_0,\s)&(\pi_1,\s)&(\pi_2,\s)\\ \hline
 &&&&&&&\\
\widetilde{\theta}_1+\widetilde{\theta}_5&2\theta(q-1)&0&-2\theta&0&0&2&-2\\
\end{array}$$\normalsize
\end{lemme}

\subsection{The outer values of~$\widetilde{\theta}_1$ and~$\widetilde{\theta}_5$}
Let~$U_0=\langle x_a\rangle\langle x_b\rangle X_{a+b}
X_{2a+b}\subseteq U$. Since~$U_0$ is~$\s$-stable, we have~$\Uz=U_0\rtimes\langle \s\rangle\subseteq\B$. 
We recall that the irreducible characters of~$B$ are given on p.~87 of~\cite{Eno}. We use the same notation. Let~$\theta_3(1)$ be the irreducible character of~$B$ of degree~$\frac{1}{2}q(q-1)^2$. Using the character tables of~$\Gf$ and~$B$ we deduce that~$\Res_B^{\Gf}\theta_5=\theta_3(1)$. Thus~$\Res_{\B}^{\G}\widetilde{\theta}_5=\widetilde{\theta}_3(1)$. We now construct~$\widetilde{\theta}_3(1)$ using~$\Uz$.
We set~$\lambda:\F_q\longrightarrow\{-1,1\}$, such that~$\lambda(x)=1$ if~$X^2+X+x$ has a root in~$\F_q$ and~$\lambda(x)=-1$ otherwise. Let~$k,l\in\{0,1\}$ and define the linear character~$\lambda(k,l)$ of~$U_0$ on the generators by~$\lambda(k,l)(x_a)= (-1)^k$, $\lambda(k,l)(x_b)=(-1)^l$ and~$\lambda(k,l)(x_{a+b}(u)x_{2a+b}(v))= \lambda(u+v)$.

\begin{lemme}
The characters~$\lambda(0,0)$ and~$\lambda(1,1)$ are $\s$-stable.
\end{lemme}

\begin{preuve}{}
Since~$\lambda(v)=1$ if and only if~$\lambda(v^{\theta})=1$ and since~$\lambda(v)=-1$ if and only if~$\lambda(v^{\theta})=-1$, it follows that $$\lambda(k,l)^{\s}(x_{a+b}(u)x_{2a+b}(v))=\lambda(k,l)(x_{a+b}(u)x_{2a+b}(v)).$$ We remark that if~$(k,l)\in\{(0,0),(1,1)\}$, then the map has the same values on~$x_a$ and~$x_b$. Thus~$\lambda(0,0)$
and~$\lambda(1,1)$ are $\s$-stable.

\end{preuve}
Using~\cite{brunat1} Lemma~3.5 the characters~$\lambda(0,0)$
and~$\lambda(1,1)$ can be extend to~$\Uz$ and we thus obtain~$\widetilde{\lambda}(0,0)$ and~$\widetilde{\lambda}(1,1)$. 

\begin{proposition}\label{conU} The conjugacy classes of~$U_0$ (resp. of~$\Uz$) are given in Table~\ref{tabconjU0} (resp.~in Table~\ref{tabconjU0ext}). In tables~\ref{tabUzInd} and~\ref{tabUzIndext}, we give the induction formula of~$\Uz$ to~$\B$.
\end{proposition}
\begin{table}
$$\begin{array}{ll|l|l}
&\text{Representatives}&\text{Number}&|\Cen_{U_0}(x)|\\ \hline
&&\\
&1&1&4q^2\\
&x_a&1&2q^2\\&x_b&1&2q^2\\&x_{a+b}(u)&q-1&4q^2\\&x_{2a+b}(v)&q-1&4q^2\\&x_ax_b&1&2q^2\\
&x_{a+b}(u)x_{2a+b}(v)&(q-1)^2&4q^2\\u\neq
0&\hspace{3pt}x_ax_{a+b}(u)&q-1&2q^2\\ u\neq
0&\hspace{3pt}x_bx_{a+b}(u)&q-1&2q^2\\ v\neq
0,1&\hspace{3pt}x_ax_{2a+b}(v)&q-2&2q^2\\v\neq
0,1&\hspace{3pt}x_bx_{2a+b}(v)&q-2&2q^2\\u\neq
0&\hspace{3pt}x_ax_bx_{a+b}(u)&q-1&2q^2\\v\neq
0,1&\hspace{3pt}x_ax_bx_{a+b}(v)&q-2&2q^2\\u,v\neq
0,1&\hspace{3pt}x_ax_{a+b}(u)bx_{2a+b}(v)&\frac{1}{2}(q-2)^2&2q^2\\ u,v\neq
0,1&\hspace{3pt}x_bx_{a+b}(u)bx_{2a+b}(v)&\frac{1}{2}(q-2)^2&2q^2\\u,v\neq
0,1&\hspace{3pt}x_ax_bx_{a+b}(u)bx_{2a+b}(v)&\frac{1}{2}(q-2)^2&2q^2
\end{array}$$
\caption{ Conjugacy classes of~$U_0$. \label{tabconjU0}}
\end{table}
\begin{table}
$$\begin{array}{ll|l|l}
 &\text{Representatives}&\text{Number}&|\Cen_{\Uz}(x)|\\ \hline
 &&\\
&(1,\s)&1&4q\\ &(x_ax_{a+b}(u),\s)&q&4q\\u\neq
0\hspace{3pt}&(x_{a+b}(u),\s)&q-1&4q
\end{array}$$
\caption{Outer classes of~$\Uz$.\label{tabconjU0ext}}
\end{table}
\begin{table}
\small
$$\renewcommand{\arraystretch}{1.4}\begin{array}{l|l}
&\Ind_{U_0}^B\phi\\ \hline
 &\\ A_1&\displaystyle{\frac{1}{4}q^2(q-1)^2\phi(1)}\\&\\
A_2&\displaystyle{\frac{1}{4}q^2(q-1)\sum_{u\in\K^{\times}}\phi(x_{2a+b}(u))}\\&\\
A_{31}&\displaystyle{\frac{1}{4}q^2(q-1)\sum_{u\in\K^{\times}}x_{a+b}(u)}\\&\\A_{32}
&\displaystyle{\frac{q^2}{4}\sum_{u,v\in\K^{\times}}\phi(x_{a+b}(u)x_{2a+b}(v))}\\&\\
A_{41}&\displaystyle{\frac{1}{2}q(q-1)+\sum_{u\neq
u+1}\phi(x_bx_{a+b}(u)x_{2a+b}(u^2)}\\&\\A_{42}&\displaystyle{\frac{q}{2}\sum_{u\in\K^{\times}}\phi(x_bx_{a+b}(u))
+\sum_{u\neq 0,1}\phi(x_bx_{2a+b}(u))+\sum_{u\neq
v^2}\phi(x_bx_{a+b}(u)x_{2a+b}(v)})\\&\\A_{51}&\displaystyle{\frac{1}{2}q(q-1)\sum_{u\neq
u+1}\phi(x_ax_{a+b}(u)x_{2a+b}(u)})\\&\\A_{52}
&\displaystyle{\frac{q}{2}\sum_{u\in\K^{\times}}\phi(x_ax_{a+b}(u))+\sum_{u\neq
0,1}\phi(x_ax_{2a+b}(u))+\sum_{u\neq
v}\phi(x_bx_{a+b}(u)x_{2a+b}(v))}\\&\\A_{61}
&\displaystyle{\sum_{\lambda(u)=1}\phi(x_ax_bx_{a+b}(u))+\sum_{\lambda(u)=1,u\neq
0}\phi(x_ax_bx_{2a+b}(u))}\\&\displaystyle{+\sum_{\lambda(u+v)=1}\phi(x_ax_bx_{a+b}(u)x_{2a+b}(v))}\\A_{62}
&\displaystyle{\sum_{\lambda(u)=-1}\phi(x_ax_bx_{a+b}(u))+\sum_{\lambda(u)=-1,u\neq
1}\phi(x_ax_bx_{2a+b}(u))}\\&\displaystyle{+
\sum_{\lambda(u+v)=-1}\phi(x_ax_bx_{a+b}(u)x_{2a+b}(v)})\\
\end{array}$$
\caption{Induction formula from~$U_0$ to~$B$\label{tabUzInd}}
\end{table}
\normalsize
\begin{table}
\small$$\begin{array}{c|c|c|c|c}&(1,\s)&(x_a,\s)&(x_{a+b},\s)&(x_ax_{a+b},\s)\\
\hline &&&& \\ \Ind_{\Uz}^{\B}\phi
 &\displaystyle{\frac{1}{2}q(q-1)\phi(1,\s)}
 &\displaystyle{\sum_{\lambda(u)=1}\phi(x_ax_{a+b}(u),\s)}
 &\displaystyle{\frac{q}{2}\sum_{u\neq 0}\phi(x_{a+b}(u),\s)}
 &\displaystyle{\sum_{\lambda(u)=-1}\phi(x_ax_{a+b}(u),\s)}
\end{array}$$
\normalsize
\caption{Outer values of the induction formula.\label{tabUzIndext}}
\end{table}

\begin{preuve}{}
We remark that~$\Ker\lambda=\{t+t^{\theta}\ |\
t\in\F_q\}$. Thus using~\cite{Eno} and the Chevalley relations, we obtain the fusion of classes. The result then follows.

\end{preuve}
\begin{lemme}
The values of~$\Ind_{\Uz}^{\B}\widetilde{\lambda}(1,1)$ are given in Table~\ref{Indlamb}.
\end{lemme}
\begin{table}
\small$$
\renewcommand{\arraystretch}{1.4}\begin{array}{l|l|l|l|l|l}
&A_1&A_2&A_{31}&A_{32}&A_{41}\\ \hline &&&&&\\
\Ind_{U_0}^B\lambda(1,1) &\displaystyle{\frac{1}{4}q^2(q-1)^2}
&\displaystyle{-\frac{1}{4}q^2(q-1)}
&\displaystyle{-\frac{1}{4}q^2(q-1)} &\displaystyle{\frac{q^2}{4}}
&\displaystyle{-\frac{1}{4}q^2(q-1)}\\ \multicolumn{6}{c}{}\\
\multicolumn{6}{c}{}\\&A_{42}&A_{51}&A_{52}&A_{61}&A_{62}\\
\cline{2-6} &&&&&\\
 &\displaystyle{\frac{q^2}{4}}
&\displaystyle{-\frac{1}{4}q^2(q-1)} &\displaystyle{\frac{q^2}{4}}
&\displaystyle{\frac{q^2}{4}} &\displaystyle{\frac{-q^2}{4}}
%\end{array}$$
%\begin{array}{l|l|l|l|l}
\\ \multicolumn{6}{c}{}\\ \multicolumn{6}{c}{}\\\multicolumn{6}{c}{}\\
&(1,\s)&(x_a,\s)&(x_{a+b},\s)&(x_ax_{a+b},\s)\\ \cline{1-5} &&&&\\
\Ind_{\Uz}^{\B}\lambda(1,1)&\displaystyle{\frac{1}{2}q(q-1)}&\displaystyle{-\frac{q}{2}}&
\displaystyle{-\frac{q}{2}}&\displaystyle{\frac{q}{2}}
\end{array}$$
\normalsize
\caption{Values of~$\Ind_{U_0}^B\lambda(1,1)$ and outer values of~$\Ind_{\Uz}^{\B}\widetilde{\lambda}(1,1)$\label{Indlamb}}
\end{table}

\begin{preuve}{}
It suffices to use Proposition~\ref{conU} and the relation~$\displaystyle{\sum_{t\in \F_q}\lambda(t)=0}$.

\end{preuve}

\begin{lemme}\label{calc}There exists an extension~$\xi$ of~$\theta_3(1)$ such that~$$\Ind_{\Uz}^{\B}\widetilde{\lambda}(1,1)=\frac{1}{4}({q+2\theta})\xi+\frac{1}{4}(q-2\theta)\xi\varepsilon.$$
\end{lemme}
\medskip
\begin{preuve}{}
We set~$\chi=\Ind_{\Uz}^{\B}\widetilde{\lambda}(1,1)$. Using the character table of~$B$ (see~\cite{Eno}~p.87) we prove that~$\langle \Res_B^{\B}\chi,\theta_3(1)\rangle_B=\frac{1}{2}q$. Thus there exist an extension~$\xi$ of~$\theta_3(1)$, integers~$n$ and~$n_\varepsilon$ where~$n\geq n_\varepsilon$, such that~$n+n_{\varepsilon}=\frac{1}{2}q$ and~$\chi=n\xi+n_\varepsilon\xi\varepsilon$.
We have~$\chi\varepsilon=n_\varepsilon\xi+n\xi\varepsilon$.
We deduce that~$\chi-\chi\varepsilon=(n-n_\varepsilon)(\xi-\xi\varepsilon).$
It follows that~$\langle
\chi-\chi\varepsilon,\chi-\chi\varepsilon\rangle_{\B}=2(n-n_\varepsilon)^2.$
Furthermore, since~$\langle
\chi-\chi\varepsilon,\chi-\chi\varepsilon\rangle_{\B}=q$ we get~$2(n-n_\varepsilon)^2=2\theta^2$ and hence~$(n-n_\varepsilon)^2=\theta^2$. This yields~$n-n_\varepsilon=\theta$.
Solving the system
$$\left\{\begin{array}{ccc} n+n_\varepsilon&=&\frac{1}{2}q\\
n-n_\varepsilon&=&\theta
\end{array}\right.$$ we find~$n=\frac{1}{4}(q+2\theta)$ and~$n_\varepsilon=\frac{1}{4}(q-2\theta)$.

\end{preuve}
\begin{lemme}\label{val1}
The outer values of~$\widetilde{\theta}_3(1)$ are: $$\begin{array}{c|c|c|c|c|c}
 &(1,\s)&(x_a,\s)&(x_{a+b},\s)&(x_ax_{a+b},\s)&(\pi_0,\s)\\ \hline
 &&&&&\\
\widetilde{\theta}_3(1)&\theta(q-1)&-\theta&-\theta&\theta&0\\
\end{array}$$
\end{lemme}
\begin{preuve}{}
Using Lemma~\ref{calc}, we compute the values of~$\xi$. For example to compute $\xi(1,\s)$ we set~$\xi(1,\s)=\alpha$. Then we have:
$$\frac{1}{4}(q+2\theta)\alpha-\frac{1}{4}({q-2\theta})\alpha=\frac{1}{2}q(q-1),$$ and we deduce~$\displaystyle{\alpha=\theta(q-1)}$. Moreover, since~$\xi(1,\s)>0$, it follows that~$\widetilde{\theta}_3(1)=\xi$.

\end{preuve}
We now set~$\phi=\Res_{\B}^{\G}(\widetilde{\theta}_1+\widetilde{\theta}_5)$. We have~$\langle \,\phi,\,\widetilde{\theta}_3(1)\,\rangle_{\B} = 1$ and~$\langle \,\phi,\,\varepsilon\widetilde{\theta}_3(1)\,\rangle_{\B}
= 0$. Moreover~$\Res_{\B}^{\G}\widetilde{\theta}_5$ is an extension of~$\theta_3(1)$. Thus~$\Res_{\B}^{\G}\widetilde{\theta}_5=\widetilde{\theta}_3(1)$ and~$\Res_{\B}^{\G}\widetilde{\theta}_1=\phi-\widetilde{\theta_3}(1)$.
We then obtain the values of~$\widetilde{\theta}_5$ and~$\widetilde{\theta_1}$
on~$(\pi_0,\s)$,~$(1,\s)$,~$(x_a,\s)$,~$(x_{a+b},\s)$
and~$(x_ax_{a+b},\s)$. This leaves us with the values of~$\widetilde{\theta}_1$ and $\widetilde{\theta}_5$
on~$(\pi_1,\s)$ and~$(\pi_2,\s)$ which still need to be computed.
\begin{lemme}\label{val2}
We have~$\widetilde{\theta}_1(\pi_1,\s)=\widetilde{\theta}_5(\pi_1,\s)=1$ and~$\widetilde{\theta}_1(\pi_2,\s)=\widetilde{\theta}_5(\pi_2,\s)=-1$.
\end{lemme}
\begin{preuve}{}
Setting~$\alpha=\widetilde{\theta}_1(\pi_1,\s)$ and~$\beta=\widetilde{\theta}_5(\pi_1,\s)$ we have~$\alpha+\beta=2$. The orthogonality relations of the rows give~$|\alpha|^2+|\beta|^2=2$. Since~$(\pi_1,\s)$ and~$(\pi_1,\s)^{-1}$ are conjugate, it follows that~$\alpha$ and~$\beta$ are real numbers. Substituting~$\beta$ by~$2-\alpha$ we see that~$\alpha$ is a root of~$X^2-2X+1$ and therefore~$\alpha=1$. It follows that~$\beta=1$. Similarly we have~$\widetilde{\theta}_1(\pi_2,\s)=\widetilde{\theta}_5(\pi_2,\s)=-1$.

\end{preuve}

\noindent In particular through Lemmas~\ref{val1} and~\ref{val2} and Propositions~\ref{chii} and~\ref{theta1} we obtain the character table of~$\G$. Thus Theorem~\ref{mainth} is proved. 
\begin{sidewaystable}
%\small
%\begin{changemargin}{-1.5cm}{0.5cm}
$$\renewcommand{\arraystretch}{1.4}
\begin{array}{lr|c|c|c|c|c|c|c}
&&(1,\s)&(x_a,\s)&(x_{a+b},\s)&(x_ax_{a+b},\s)&(\pi_0,\s)&(\pi_1,\s)&(\pi_2,\s)\\\cline{3-9}
|C_{\G}|&&2q^2(q-1)(q^2+1)&4q&2q^2&4q&2(q-1)&2(q+2\theta+1)&2(q-2\theta+1)\\\hline
1&1&1&1&1&1&1&1&1\\
%&&&&&&&\\
\widetilde{\theta}_4&1&q^2&0&0&0&1&-1&-1\\
%&&&&&&&\\
\widetilde{\theta}_1&1&\theta(q-1)&\theta&-\theta&-\theta&0&1&-1\\
%&&&&&&&\\
\widetilde{\theta}_5&1&\theta(q-1)&-\theta&-\theta&\theta&0&1&-1\\
%&&&&&&&\\
\widetilde{\chi}_{\pi_0}(i)&i\in
E_0&q^2+1&1&1&1&\varepsilon_0^{i}(\pi_0)&0&0\\
%&&&&&&&\\
\widetilde{\chi}_{\pi_1}(j)&j\in
E_1&(q-1)(q-2\theta+1)&-1&2\theta-1&-1&0&-\varepsilon_1^j(\pi_1)&0\\
%&&&&&&&\\
\widetilde{\chi}_{\pi_2}(k)&k\in
E_2&(q-1)(q+2\theta+1)&-1&-2\theta-1&-1&0&0&-\varepsilon_2^k(\pi_2)\\
\end{array}$$
%\end{changemargin}
\caption{Values of the outer characters~$B_2(q)\semi{\s}$ on outer classes\label{CaractereB2}}
\end{sidewaystable}
\normalsize
\section{Shintani Descents}
In this section we use the results and notation of~\S\ref{DL} and~\S\ref{B2}. Just recall that~$\Galg$ is a simple group of type~$B_2$ and that~$F$ is the generalized Frobenius map such that~$\Galg^F=\Su(q)$~(where~$q=2^{2n+1}$). Recall furthermore that~$\delta=2$. We will first give Shintani descents between~$\G$ and~$\Su(q)$ and then obtain some results on the unipotent characters of~$\Su(q)$.
\subsection{Shintani correspondence}
\begin{proposition}\label{parite} We have:
\begin{itemize}
\item[\textbullet] If~$n$ is odd, then~$N_{F/F^2}\Cl(x_ax_bx_{a+b})=\Cl(x_a,\s)$.
\item[\textbullet] If~$n$ is even, then~$N_{F/F^2}\Cl(x_ax_bx_{a+b})=\Cl(x_ax_{a+b},\s)$.
\end{itemize}
\end{proposition}
\begin{preuve}{}
We search for a~$g\in\Galg$ such that~$x_a=g^{-1}F(g)$. To this end suppose there exist~$u,v,w,t\in\overline{\F}_2$ such that~$g=x_a(u)x_b(v)x_{a+b}(w)x_{2a+b}(t)$.
Then using the Chevalley relations of~$\Galg$ we immediately deduce:
$$\begin{array}{lcl}
%g^{-1}&=&x_a(u)x_b(v)x_{a+b}(uv+w)x_{2a+b}(u^2v+t),\\
F(g)&=&x_a(v^{2^n})x_b(u^{2^{n+1}})x_{a+b}(t^{2^n}+v^{2^n}u^{2^{n+1}})x_{2a+b}(w^{2^{n+1}}+v^{2^{n+1}}u^{2^{n+1}}),
\\gx_a&=&x_a(u+1)x_b(v)x_{a+b}(w+v)x_{2a+b}(t+v).
\end{array}$$ By uniqueness of this decomposition we obtain:
$$\left\{\begin{array}{lcl} u+1&=&v^{2^n}\\ v&=&u^{2^{n+1}}\\
t^{2^n}+v^{2^{n}}u^{2^{n+1}}+w+v&=&0\\w^{2^{n+1}}+v^{2^{n+1}}u^{2^{n+1}}+t+v&=&0\\
\end{array}\right.$$
that is:$$\left\{\begin{array}{lcl} u+1&=&v^{2^n}\\
v&=&u^{2^{n+1}}\\ t^{2^n}+w+uv&=&0\\w^{2^{n+1}}+t+u^2v&=&0.\\
\end{array}\right.$$
Suppose this system has a solution~$(u,v,w,t)$. Then we have:
$$\left\{\begin{array}{lcl} u^q+u+1&=&0\\ v^q+v+1&=&0\\
w^q+w+v&=&0\\t^q+t+v&=&0.\\
\end{array}\right.$$
Using these relations we find~$gF^2(g^{-1})=x_ax_bx_{a+b}(u+v+1)x_{2a+b}(u^2+v+1)$.
We now prove that this system has a solution. We write~$h\in\overline{\F}_2$ for a root of~$X^2+X+1$ and~$k\in\overline{\F}_2$ for a root of~$X^2+X+h^2$ and we prove by induction on~$m$ that:
$$\left\{\begin{array}{lclll} h^{2^m}&=&h+1&&\text{if}\ m\ \text{odd}\\ h^{2^m}&=&h&&\text{if}\ m\ \text{is even}
\end{array}\right.$$
and: $$\left\{\begin{array}{lclll} k^{2^m}&=&k&&\text{if}\ m =0\mod
4\\ k^{2^m}&=&k+h^2&&\text{if}\ m =1\mod 4\\
k^{2^m}&=&k+1&&\text{if}\ m =2\mod 4\\ k^{2^m}&=&k+h&&\text{if}\ m
=3\mod 4.\end{array}\right.$$ We now study each of these cases:
\begin{itemize}
\item Suppose~$n=1\mod 4$. In this case~$n$ is odd. We set~$v=u=h$ and~$t=w=k$. We get~$t^{2^n}+t+uv=0$ and~$t^{2^{n+1}}+t+u^2v=0$.
The element~$g=x_a(h)x_b(h)x_{a+b}(k)x_{2a+b}(k)$ is a solution and we have~$gF^2(g^{-1})=x_ax_bx_{a+b}$.
\item
Suppose~$n=3\mod 4$. In this case~$n$ is odd. We set~$v=u=h$ and~$t=w=h+k$. We obtain~$t^{2^n}+w+uv=0
$ and~$w^{2^{n+1}}+t+u^2v=0$.
The element~$g=x_a(h)x_b(h)x_{a+b}(h+k)x_{2a+b}(h+k)$ is a solution and~$gF^2(g^{-1})=x_ax_bx_{a+b}$.
\item
Suppose~$n=0\mod 4$. In this case~$n$
is even. We set~$u=h$,~$v=u+1$,~$t=k$ and~$w=k+1$.  The element $g=x_a(h)x_b(h+1)x_{a+b}(k+1)x_{2a+b}(k)$
is a solution and we have~$gF^2(g^{-1})=x_ax_bx_{2a+b}$.
\item
Suppose~$n=2\mod 4$. In this case~$n$
is even. We put~$u=h$,~$v=u+1$,~$t=k$ and~$w=k$.  Then the element~
$g=x_a(h)x_b(h+1)x_{a+b}(k)x_{2a+b}(k)$ is a solution and we have~$gF^2(g^{-1})=x_ax_bx_{2a+b}$.
\end{itemize}
Using the definition of the Shintani correspondence (see~\S\ref{secshin}), the claim is proved.

\end{preuve}

\subsection{Shintani descents of unipotent characters}
We fix~$i$ a primitive fourth root of unity. We recall that~$\rho_0=x_{a}x_{a+b}$ and that~$\St$ is the Steinberg character of~$\Su(q)$. We denote by~$\mathcal{W}$ the cuspidal unipotent character of degree~$\theta(q-1)$ such that~$\mathcal{W}(\rho_0)=\theta i$. We write~$1_{\G}$,~$\widetilde{\theta}_1$,~$\widetilde{\theta}_4$ and~$\widetilde{\theta}_5$ for the extensions of the unipotent characters of~$\Gf$ as above. 
\begin{theoreme}\label{Shintani} We have~$\Sh_{F^2/F} 1_{\G} = 1_{\Su(q)}$ and~$Sh_{F^2/F}\widetilde{\theta}_4 = \St$. Setting~$\zeta_0=\frac{\sqrt{2}}{2}(-1-i)$ we have:
\begin{itemize}
\item[\textbullet] If $n$ is even, then: $$\begin{array}{lll}
\Sh_{F^2/F}\widetilde{\theta}_5&=&-\frac{\zeta_0}{\sqrt{2}}
\mathcal{W}-\frac{\overline{\zeta_0}}{\sqrt{2}}\overline{\mathcal{W}}\\
\Sh_{F^2/F}\widetilde{\theta}_1&=&-\frac{\overline{\zeta_0}}{\sqrt{2}}\mathcal{W}-\frac{\zeta_0}{\sqrt{2}}\overline{\mathcal{W}}
\end{array}$$
\item[\textbullet] If~$n$ is odd, then: $$\begin{array}{lll}
\Sh_{F^2/F}\widetilde{\theta}_1&=&-\frac{\zeta_0}{\sqrt{2}}
\mathcal{W}-\frac{\overline{\zeta_0}}{\sqrt{2}}\overline{\mathcal{W}}\\
\Sh_{F^2/F}\widetilde{\theta}_5&=&-\frac{\overline{\zeta_0}}{\sqrt{2}}\mathcal{W}-\frac{\zeta_0}{\sqrt{2}}\overline{\mathcal{W}}
\end{array}$$
\end{itemize}
\end{theoreme}

\begin{preuve}{}
We say that a class of~$\Su(q)$ is of type~$\pi_i$ ($i\in\{0,1,2\}$) if
there exists some~$k\in E_i$ such that~$\pi_i^k$ is a representative   
of the class. Similarly we say that a class of~$\G$ is of              
type~$\pi_i$ ($i\in\{0,1,2\}$) if there exists some~$k\in E_i$          
such that~$(\pi_i^k,\s)$ is a representative of the class.              
Using Relation~\ref{eqcardinal} in~\S\ref{secshin} and Theorem~\ref{tconj}, we            
see that~the classes of~$\Su(q)$ of type~$\pi_0$,~$\pi_1$             
and~$\pi_2$ are sent by~$N_{F/F^2}$ to classes of~$\G$ of            
type~$(\pi_0,\s)$,~$(\pi_1,\s)$ and~$(\pi_2,\s)$ respectively.         
Since~$1_{\G}$,~$\widetilde{\theta}_1$,~$\widetilde{\theta}_4$         
and~$\widetilde{\theta}_5$ are constant on these classes, we do not    
explicitly need to know the correspondence of these classes to compute the     
Shintani descents of these characters. Using a similar argument we      
obtain that~$N_{F/F^2}\Cl(1)=\Cl(1,\s)$ and~$N_{F/F^2}\Cl(x_{a+b}      
x_{2a+b})=\Cl(x_{a+b},\s)$. We now use Proposition~\ref{parite}        
\begin{itemize}
\item[\textbullet] Suppose that~$n$ is odd. Then we get $$N_{F/F^2}\Cl(\rho_0)=\Cl(x_a,\s)\quad\textrm{and}\quad N_{F/F^2}\Cl(\rho_0^{-1})=\Cl(x_ax_{a+b},\s).$$ Using the values of outer characters obtained in Theorem~\ref{mainth} we have:
%\small
$$\renewcommand{\arraystretch}{1.4}
\begin{array}{l|c|c|c|c|c|c|c}
&1&\rho_0&\s_0&\rho_0^{-1}&\pi_0^i&\pi_1^j&\pi_2^k\\\hline
\Sh_{F^2/F}1_{\G}&1&1&1&1&1&1&1\\%&&&&&&&\\
\Sh_{F^2/F}\widetilde{\theta}_4&q^2&0&0&0&1&-1&-1\\
%&&&&&&&\\
\Sh_{F^2/F}\widetilde{\theta}_1&\theta(q-1)&\theta&-\theta&-\theta&0&1&-1\\
%&&&&&&&\\
\Sh_{F^2/F}\widetilde{\theta}_5&\theta(q-1)&-\theta&-\theta&\theta&0&1&-1\\
\end{array}$$
%\caption{Caract\`eres Ext\'erieurs de $\G$ sur les classes
%ext\'erieures.\label{table}}
%\end{table}
%\normalsize
This yields~$\Sh_{F^2/F}1_{\G} = 1$ and~$\Sh_{F^2/F}\widetilde{\theta_4}=\St$.
% Le calcul des produits
%scalaires avec les caractères irréductibles de~$\Su$, montre que
%les constituants de~$Sh_{F^2/F}\widetilde{\theta_1}$ et de
%$Sh_{F^2/F}\widetilde{\theta_5}$ sont des caractères unipotents de
%$\Su$, plus précisement, on a
Moreover, using the character Table of~$\Su(q)$ we compute:
$$
\begin{array}{ll}
\langle\,\Sh_{F^2/F}\widetilde{\theta}_1,1\,\rangle_{\Su}=0&\langle\,\Sh_{F^2/F}\widetilde{\theta}_5,1\,\rangle_{\Su}=0\\
\langle\,\Sh_{F^2/F}\widetilde{\theta}_1,\St\,\rangle_{\Su}=0&\langle\,\Sh_{F^2/F}\widetilde{\theta}_5,\St\,\rangle_{\Su}=0\\
\langle\,\Sh_{F^2/F}\widetilde{\theta}_1,\mathcal{W}\,\rangle_{\Su}=-\zeta_0\sqrt{2}/2&\langle\,\Sh_{F^2/F}\widetilde{\theta}_5,\mathcal{W}\,\rangle_{\Su}=-\overline{\zeta_0}\sqrt{2}/2\\
\langle\,\Sh_{F^2/F}\widetilde{\theta}_1,\overline{\mathcal{W}}\,\rangle_{\Su}=-\overline{\zeta_0}\sqrt{2}/2&\langle\,\Sh_{F^2/F}\widetilde{\theta}_5,\overline{\mathcal{W}}\,\rangle_{\Su}=-\zeta_0\sqrt{2}/2\\
\end{array}
 $$ We thus deduce that: $$\begin{array}{lll}
\Sh_{F^2/F}\widetilde{\theta}_1&=&-\zeta_0 \sqrt{2}/2
\mathcal{W}-\overline{\zeta_0}\sqrt{2}/2\overline{\mathcal{W}}\\
\Sh_{F^2/F}\widetilde{\theta}_5&=&-\overline{\zeta_0}\sqrt{2}/2\mathcal{W}-\zeta_0\sqrt{2}/2\overline{\mathcal{W}}\\
\end{array}$$
\item[\textbullet] If~$n$ is even, we proceed similarly using the identities~$N_{F/F^2}\Cl(\rho_0)=\Cl(x_ax_{a+b},\s)$ and~$N_{F/F^2}\Cl(\rho_0^{-1})=\Cl(x_a,\s)$.
\end{itemize}

\end{preuve}
\subsection{Roots associated to the unipotent characters of~$\Su(q)$}
We denote by~$\omega_{\mathcal{W}}$ and~$\omega_{\overline{\mathcal{W}}}$ the roots of unity associated to~$\mathcal{W}$ and~$\overline{\mathcal{W}}$ as in~\S\ref{root}.
In~\cite{Lusztig}~\S$7.4$ G.~Lusztig shows that~$
\{\omega_{\mathcal{W}},\omega_{\overline{\mathcal{W}}}\}=\{\zeta_0,\overline{\zeta}_0\}$. We now make this result more precise.

\begin{theoreme}\label{RacineB2}
The roots associated to~$\mathcal{W}$ and~$\overline{\mathcal{W}}$ are:
$$\renewcommand{\arraystretch}{1.4}\begin{array}{l|ll}
&\mathcal{W}&\overline{\mathcal{W}}\\\hline n\ \textrm{
odd}&\zeta_0&\overline{\zeta_0}\\ n\ \textrm{
even}&\overline{\zeta_0}&\zeta_0\\
\end{array}$$
\end{theoreme}

\begin{preuve}{} To compute the almost characters of~$\Su(q)$ we need to know the generalized Deligne-Lusztig characters~${R}_{w}$. The $F$-classes of~$W$ have the representatives~$1$,~$w_a$ and~$w_aw_bw_a$.
The Suzuki group with parameter~$q$ therefore has~three Deligne-Lusztig characters with the degrees:
$$\begin{array} {lll} R_1(1)=q^2+1,& R_{w_a}(1)=
(q-1)(q-r+1),&R_{w_aw_bw_a}(1)=(q-1)(q+r+1).
\end{array}$$
These characters are explicitly computed in~\cite{Geck}. We recall that (see~\cite{Geck}~Prop.~$4.6.7$):
$$\begin{array}{lll}
R_{1}&=&1+\St\\ R_{w_a}&=&1-\mathcal{W}-\overline{\mathcal{W}}-\St\\
R_{w_aw_bw_a}&=&1+\mathcal{W}+\overline{\mathcal{W}}-\St.
\end{array}$$
The Weyl group~$W$ has three~$F$-stable characters denoted by~$\rho_1$,~$\rho_2$ and~$\rho_3$. Let~$\widetilde{\rho}_1$,~$\widetilde{\rho}_2$ and~$\widetilde{\rho}_3$ be their extensions to~$W\semi{F}$ such that:
$$\begin{array}{c|ccc}
&(1,F)&(w_a,F)&(w_aw_bw_a,F)\\\hline
\widetilde{\rho}_1&1&1&1\\
\widetilde{\rho}_2&1&-1&-1\\
\widetilde{\rho}_3&0&-\sqrt{2}&\sqrt{2}\\
\end{array}$$ The three almost characters of~$\Su(q)$ corresponding to these extensions are~(\S\ref{FourierMat}):
$$\begin{array}{lll} \RR_{\widetilde{\rho}_1}&=&1_{\Su(q)}\\ \RR_{\widetilde{\rho}_2}&=&\St\\
\RR_{\widetilde{\rho}_3}&=&\frac{\sqrt{2}}{2}(\mathcal{W}+\overline{\mathcal{W}}).\\
\end{array}$$
Therefore we have:
$$\begin{array}{lll}\langle\,\RR_{\widetilde{\rho}_3},\mathcal{W}\,\rangle_{\Su}=\sqrt{2}/2&\textrm{and}&
\langle\,\RR_{\widetilde{\rho}_3},\overline{\mathcal{W}}\,\rangle_{\Su}=\sqrt{2}/2.\end{array}$$
Moreover~$\theta_1$ is in the principal series of~$\Gf$. Using Theorem~\ref{thDM} we deduce that:
$$\sqrt{2}\cyc{\Sh_{F^2/F}\widetilde{\theta}_1,\mathcal{W}}_{\Su}=\pm \omega_{\mathcal{W}}.$$
The sign is due to the fact that~$\Sh_{F^2/F}\widetilde{\theta}_1=\pm\Sh_{F^2/F}\chi_{\widetilde{\rho}_3}$. Since~$\omega_{\mathcal{W}}$ is either~$\zeta_0$ or~$\overline{\zeta}_0$ and using that~$\overline{\zeta}_0\neq -\zeta_0$, we can obtain the root. Indeed
\begin{itemize}
\item Either~$\sqrt{2}\cyc{\Sh_{F^2/F}\widetilde{\theta}_1,\mathcal{W}}$ is~$\zeta_0$ or~$\overline{\zeta}_0$ and in this case~$\omega_{\mathcal{W}}=\sqrt{2}\cyc{\Sh_{F^2/F}\widetilde{\theta}_1,\mathcal{W}}$.
\item Or~$\sqrt{2}\cyc{\Sh_{F^2/F}\widetilde{\theta}_1,\mathcal{W}}$ is neither~$\zeta_0$ nor~$\overline{\zeta}_0$ and then~$\omega_{\mathcal{W}}=-\sqrt{2}\cyc{\Sh_{F^2/F}\widetilde{\theta}_1,\mathcal{W}}$.
\end{itemize} Using~\ref{parite} the claim is proved.

\end{preuve}

\subsection{Fourier matrices}

The unipotent characters of~$\Su(q)$ are distributed in three families~$\mathcal{F}_1=\{1\}$, $\mathcal{F}_2=\{\St\}$ and~$\mathcal{F}_3=\{\mathcal{W},\mathcal{\overline{W}}\}$.

\begin{proposition}\label{FourierB2}
The Fourier matrices~$M_i$ ($i=1,2,3$) associated to the~$\mathcal{F}_i$ can be define up to a normalization by
$$M_1=M_2=[1],\quad\textrm{and}\quad M_3=\begin{bmatrix}\sqrt{2}/2&\sqrt{2}/2\\\sqrt{2}/2&-\sqrt{2}/2\end{bmatrix}.$$
\end{proposition}

\begin{preuve}{}
Theorem~\ref{Shintani} shows that the irreducible components of~$\Sh_{F^2/F}\widetilde{\theta}_1$ and of~$\Sh_{F^2/F}\widetilde{\theta}_5$ are the elements of the family~$\mathcal{F}_3$. On the other hand we deduce from Theorem~\ref{Shintani} that:
$$\begin{array}{ll} \text{If $n$ even:}&i \Sh_{F^2/F}\widetilde{\theta_5}=-\overline{\zeta_0} \sqrt{2}/2
\mathcal{W}+\zeta_0\sqrt{2}/2\overline{\mathcal{W}}\\ \\ \text{If $n$ is odd:}&i
\Sh_{F^2/F}\widetilde{\theta_5}=\zeta_0\sqrt{2}/2
\mathcal{W}-\overline{\zeta_0}\sqrt{2}/2\overline{\mathcal{W}}.
\end{array}$$
Using Conjecture~\ref{conjDM} we can define the Fourier matrix associated to~$\mathcal{F}_3$ as claimed.

\end{preuve}

\section*{Appendix}
\begin{table}
$$\renewcommand{\arraystretch}{1.5}\begin{array}{|l|l|c|l|}\hline
\textrm{Notation}&\textrm{Representatives}&\textrm{Number}&\textrm{Centralizer's
order}\\
\hline
A_1&h(1,1)&1&q^4(q^2-1)(q^4-1)\\
A_2&x_{2a+b}&1&q^4(q^2-1)\\
A_{31}&x_{a+b}&1&q^4(q^2-1)\\
A_{32}&x_{a+b}x_{2a+b}&1&q^4\\
A_{41}&x_{a}x_b&1&2q^2\\
A_{42}&x_{a}x_bx_{2a+b}&1&2q^2\\
B_1(i,j)&h(\gamma^i,\gamma^j)&\frac{1}{8}(q-2)(q-4)&(q-1)^2\\
B_2(i)&h(\tau'^i,\tau'^{qi})&\frac{1}{4}q(q-2)&q^2-1\\
B_3(i,j)&h(\gamma^i,\nu^{j})&\frac{1}{4}q(q-2)&q^2-1\\
B_4(i,j)&h(\nu^i,\nu^j)&\frac{1}{8}q(q-2)&(q+1)^2\\
B_5(i)&h(\tau^i,\tau^{qi})&\frac{1}{4}q^2&q^2+1\\
C_1(i)&h(1,\gamma^i)&\frac{1}{2}(q-2)&q(q-1)(q^2-1)\\
C_2(i)&h(\gamma^i,\gamma^{-i})&\frac{1}{2}(q-2)&q(q-1)(q^2-1)\\
C_3(i)&h(1,\nu^i)&\frac{1}{2}q&q(q+1)(q^2-1)\\
C_4(i)&h(\nu^i,\nu^{-i})&\frac{1}{2}q&q(q+1)(q^2-1)\\
C_1(i)&h(1,\gamma^i)x_{2a+b}&\frac{1}{2}(q-2)&q(q-1)\\
C_2(i)&h(\gamma^i,\gamma^{-i})x_{a+b}&\frac{1}{2}(q-2)&q(q-1)\\
C_3(i)&h(1,\nu^i)x_{2a+b}&\frac{1}{2}q&q(q+1)\\
C_4(i)&h(\nu^i,\nu^{-i})x_{a+b}&\frac{1}{2}q&q(q+1)\\
\hline
\end{array}$$
\caption{Conjugacy classes of~$B_2(q)$\label{ConjSp4}}
\end{table}

Let~$n$ be a non-negative integer and write~$q=2^{2n+1}$. The character table of~$B_2(q)$ is given in~\cite{Eno}~p.~$93$. In Table~\ref{ConjSp4} we recall the conjugacy classes of~$B_2(q)$.

\begin{remarque} The two classes of~$B_2(q)$ whose centralizers are of order~$2q^2$ are~$\Cl(x_ax_b)$ and~$\Cl(x_ax_bx_{a+b})$. Indeed suppose there exists~$\alpha\in\F_q$ such that $P(\alpha)=0$, where~$P=X^2+X+1$.
We can suppose that~$\alpha\not\in\{ 0,1\}$ because~$0$ and~$1$ are no roots of~$P$. Moreover~$1-\alpha^3=(1-\alpha)(\alpha^2+\alpha+1)=0$ because~$\alpha$
is a root of~$P$. We deduce that the order of~$\alpha$ divides~$3$. Furthermore~$\alpha\neq 1$, thus the order of~$\alpha$ is~$3$. It follows that~$3$ divides~$(q-1)$, that is~$q-1=0\mod 3$ which is false. Thus~$P$ is irreducible.

\end{remarque}
We set~$\alpha_{i}=\gamma_0^i+\gamma_0^{-i}$ and~$\beta_i=\nu_0^k+\nu_0^{-k}$.
In Table~\ref{CarSp4}, we recall the irreducible characters that we need in this work and correct some errors of~\cite{Eno}. 
\begin{sidewaystable}
%\begin{table}
\small
$$\begin{array}{c}
\renewcommand{\arraystretch}{1.5}
\begin{array}{l|c|c|c|c|c|c|c|c|c|c}
&A_1&A_2&A_{31}&A_{32}&A_{41}&A_{42}&B_1(i,j)&B_2(i)&B_{3}(i,j)&B_4(i,j)\\
\hline
\theta_1&q(q+1)^2/2&q(q+1)/2&q(q+1)/2&q/2&q/2&-q/2&2&0&0&0\\
\theta_4&q^4&0&0&0&0&0&1&-1&-1&1\\
\theta_5&q(q-1)^2/2&-q(q-1)/2&-q(q-1)/2&q/2&q/2&-q/2&0&0&0&-2\\
\chi_1(k,l)&(q+1)^2(q^2+1)&(q+1)^2&(q+1)^2&2q+1&1&1&\alpha_{ik}\alpha_{jl}+\alpha_{il}\alpha_{jk}&0&0&0\\
\chi_4(k,l)&(q-1)^2(q^2+1)&(q-1)^2&(q-1)^2&-(2q-1)&1&1&0&0&0&\beta_{ik}\beta_{jl}+\beta_{il}\beta_{jk}\\
\chi_5(k)&(q^2-1)^2&-(q^2-1)&-(q^2-1)&1&1&1&0&0&0&0\\
\end{array}
\\ \\ \\ \\ \\
\renewcommand{\arraystretch}{1.5}
\begin{array}{l|c|c|c|c|c|c|c|c|c}
&B_5(i)&C_1(i)&C_2(i)&C_3(i)&C_4(i)&D_1(i)&D_2(i)&D_3(i)&D_4(i)\\
\hline
\theta_1&-1&q+1&q+1&0&0&1&1&0&0\\
\theta_4&1&q&q&-q&-q&0&0&-1&-1\\
\theta_5&1&0&0&q-1&q-1&0&0&-1&-1\\
\chi_1(k,l)&0&(q+1)(\alpha_{ik}+\alpha_{il})&(q+1)\alpha_{ik}\alpha_{il}&0&0&\alpha_{ik}+\alpha_{il}&\alpha{ik}\alpha_{il}&0&0\\
\chi_4(k,l)&0&0&0&-(q-1)(\beta_{ik}+\beta_{il})&-(q-1)\beta_{ik}\beta_{il}&0&0&\beta_{ik}+\beta_{il}&\beta_{ik}\beta_{il}\\
\chi_5(k)&\tau^{ik}+\tau^{-ik}+\tau^{ikq}+\tau^{-ikq}&0&0&0&0&0&0&0&0\\
\end{array}
\end{array}$$
\normalsize \caption{\label{B2eno} Character table of~$B_2(q)$\label{CarSp4}}
\end{sidewaystable}

\noindent We recall the character table of~$\Su(q)$ wich is computed in~\cite{Suzuki} in Table~\ref{carSu}. Here we set~$\theta=2^n$.

\begin{sidewaystable}
\renewcommand{\arraystretch}{1.4}
$$\begin{array}{l|c|c|c|c|c|c|c}
&1&\s_0&\rho_0&\rho_0^{-1}&\pi_0^l&\pi_1^l&\pi_2^l\\
\hline
1_{\Su(q)}&1&1&1&1&1&1&1\\
\St&q^2&0&0&0&1&-1&-1\\
\mathcal{W}&\theta(q-1)&-\theta&\theta\sqrt{-1}&-\theta\sqrt{-1}&0&1&-1\\
\overline{\mathcal{W}}&\theta(q-1)&-\theta&-\theta\sqrt{-1}&\theta\sqrt{-1}&0&1&-1\\
X_i&q^2+1&1&1&1&\varepsilon_0^{i}(\pi_0^l)&0&0\\
Y_j&(q-2\theta+1)(q-1)&2\theta-1&-1&-1&0&-\varepsilon_1^{j}(\pi_1^l)&0\\
Z_k&(q+2\theta+1)(q-1)&-2\theta-1&-1&-1&0&0&-\varepsilon_2^{k}(\pi_2^l)\\
\end{array}$$
\normalsize
\caption{\label{carSu} Character table of~$\Su(q)$}
\end{sidewaystable}

\end{document}